\documentclass[11pt]{article}
\usepackage{amsmath}
\usepackage{dsfont}
\usepackage{mathrsfs}
\usepackage{amsmath,amssymb}
\usepackage{amsfonts}
\usepackage{hyperref}
\usepackage{amsthm}
\usepackage{graphicx}
\usepackage{subfigure}
\usepackage{xcolor}

\renewcommand{\qed}{\hfill\small{$\square$}\normalsize}

\hfuzz=\maxdimen
\theoremstyle{definition}
\newtheorem{lemma}{Lemma}[section]
\newtheorem{definition}[lemma]{Definition}

\newtheorem{theorem}[lemma]{Theorem}
\newtheorem{corollary}[lemma]{Corollary}

\newtheorem{remark}{Remark}

\numberwithin{equation}{section}
\renewcommand{\proof}{\textbf{Proof. }}
\renewcommand{\qed}{\hfill\small{$\square$}\normalsize}

\DeclareFixedFont{\Acknowledgment}{OT1}{cmr}{bx}{n}{14pt}
\begin{document}

\title{\bf Rigidity of inversive distance circle packings revisited}
\author{Xu Xu}
\maketitle

\begin{abstract}
Inversive distance circle packing metric was introduced by P Bowers and K Stephenson \cite{BS} as a generalization of
Thurston's circle packing metric \cite{T1}. They conjectured that the inversive distance circle packings are rigid.
For nonnegative inversive distance, Guo \cite{Guo} proved the infinitesimal rigidity and then Luo \cite{L3}
proved the global rigidity.
In this paper, based on an observation of Zhou \cite{Z},
we prove this conjecture for inversive distance in $(-1, +\infty)$
by variational principles.
We also study the global rigidity of a combinatorial curvature introduced in \cite{GJ4,GX4,GX6}
with respect to the inversive distance circle packing metrics where the inversive distance is in $(-1, +\infty)$.
\end{abstract}

\section{Introduction}
\subsection{Background}
In his work on constructing hyperbolic structure on 3-manifolds, Thurston (\cite{T1}, Chapter 13) introduced
the notion of circle packing metric on triangulated surfaces
with non-obtuse intersection angles.
The requirement of prescribed intersection angles corresponds to the fact that the intersection angle of
two circles is invariant under the M\"{o}bius transformations.
For triangulated surfaces with Thurston's circle packing metrics, there are singularities at the vertices.
The classical combinatorial curvature $K_i$ is introduced to describe the singularity at the vertex $v_i$,
which is defined as the angle deficit at $v_i$.
Thurston's work generalized Andreev's work on circle packing metrics on a sphere \cite{An1,An2} and
gave a complete characterization of the space of the classical combinatorial curvature.
As a corollary, he obtained the combinatorial-topological obstacle for the existence of a constant
curvature circle packing with non-obtuse intersection angles,
which could be written as combinatorial-topological inequalities.
Zhou \cite{Z} recently generalized Andreev-Thurston Theorem to the case that the intersection angles are in $[0,\pi)$.
Chow and Luo \cite{CL1} introduced a combinatorial Ricci flow, a combinatorial analogue of the smooth
surface Ricci flow, for triangulated surfaces with Thurston's circle packing metrics and established the equivalence between the
existence of a constant curvature circle packing metric and the convergence of the combinatorial Ricci flow.

Inversive distance circle packing on triangulated surfaces
was introduced by Bowers and Stephenson \cite{BS}
as a generalization of Thurston's circle packing.
Different from Thurston's circle packing, adjacent circles in inversive distance circle packing
are allowed to be disjoint and the relative distance of the adjacent circles is measured by the
inversive distance, which is a generalization of intersection angle.
See Bowers-Hurdal \cite{BH}, Stephenson \cite{St} and Guo \cite{Guo} for more information.
The inversive distance circle packings have practical applications in medical imaging and computer graphics, see \cite{HS, ZG, ZGZLYG} for example.
Bowers and Stephenson \cite{BS} conjectured that the inversive distance circle packings are rigid.
Guo \cite{Guo} proved the infinitesimal rigidity and then Luo \cite{L3} solved affirmably the conjecture for
nonnegative inversive distance with Euclidean and hyperbolic background geometry.
For the spherical background geometry, Ma and Schlenker \cite{MS} had a counterexample showing that there is in general no rigidity
and John C. Bowers and Philip L. Bowers \cite{BB} obtained a new construction
of their counterexample using the inversive geometry of the 2-sphere.
John Bowers, Philip Bowers and Kevin Pratt \cite{BBP} recently proved
the global rigidity of convex inversive distance circle packings in the Riemann sphere.
Ge and Jiang \cite{GJ2,GJ3} recently studied the deformation of combinatorial curvature
and found a way to search for inversive distance circle packing metrics with constant cone angles.
They also obtained some results on the image of curvature map for inversive distance circle packings.
Ge and Jiang \cite{GJ4} and Ge and the author \cite{GX6} further extended a combinatorial curvature introduced by Ge and the author in \cite{GX4,GX3,GX5}
to inversive distance circle packings and studied the rigidity and deformation of the curvature.

In this paper, based on an obversion of Zhou \cite{Z}, we prove Bowers and Stephenson's rigidity conjecture
for inversive distance in $(-1, +\infty)$. The main tools are the variational principle established by Guo \cite{Guo} for inversive distance circle
packings and the extension of locally convex function introduced by Bobenko, Pinkall and Springborn \cite{BPS} and systematically developed by
Luo \cite{L3}.
We refer to Glickenstein \cite{G} for a nice geometric interpretation of the variational principle in \cite{Guo}.
There are many other works on variational principles on circle packings. See Br\"{a}gger \cite{Br}, Rivin \cite{R},
Leibon \cite{L}, Chow-Luo \cite{CL1}, Bobenko-Springborn \cite{BS}, Marden-Rodin \cite{MR},
Spingborn \cite{Sp}, Stephenson \cite{St}, Luo \cite{L4}, Guo-Luo \cite{GL}, Dai-Gu-Luo \cite{DGL}, Guo\cite{Guo0} and others.

\subsection{Inversive distance circle packings}
In this subsection, we briefly recall the inversive distance circle packing introduced by Bowers and Stephenson \cite{BS}
in Euclidean and hyperbolic background geometry.
For more information on inversive distance circle packing metrics, the readers can refer to  Stephenson \cite{St},
Bowers and Hurdal \cite{BH} and Guo \cite{Guo}.

Suppose $M$ is a closed surface with a triangulation $\mathcal{T}=\{V,E,F\}$,
where $V,E,F$ represent the sets of vertices, edges and faces respectively.
Let $I: E\rightarrow (-1,+\infty)$ be a function assigning each edge $\{ij\}$ an inversive distance $I_{ij}\in (-1,+\infty)$,
which is denoted as $I>-1$ in the paper.
The triple $(M, \mathcal{T}, I)$ will be referred to as a weighted triangulation of $M$ below.
All the vertices are ordered one by one, marked by $v_1, \cdots, v_N$, where $N=|V|$
is the number of vertices, and we often use $i$ to denote the vertex $v_i$ for simplicity below.
We use $i\sim j$ to denote that the
vertices $i$ and $j$ are adjacent, i.e., there is an edge $\{ij\}\in E$ with $i$, $j$ as end points.
All functions $f: V\rightarrow \mathbb{R}$ will be regarded as column
vectors in $\mathbb{R}^N$ and $f_i=f(v_i)$ is the value of $f$ at $v_i$. And we use $C(V)$ to denote the set of functions
defined on $V$. $\mathbb{R}_{>0}$ denotes the set of positive numbers in the paper.

\begin{figure}
  \centering
  \includegraphics[width=1\textwidth]{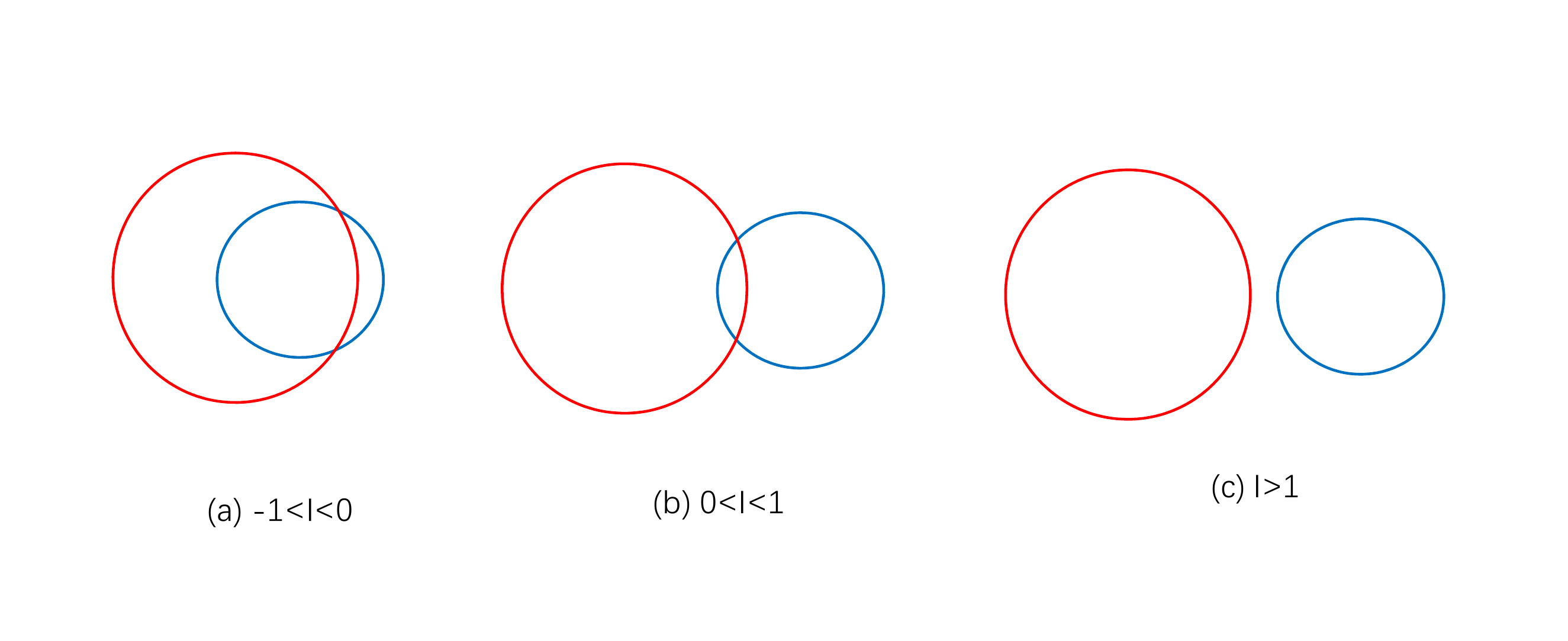}\\
  \caption{Inversive distance circle packings}\label{Inversive_distance_circle_packing}
\end{figure}

Each map $r:V\rightarrow (0,+\infty)$ is a circle packing, which could be taken as the radius $r_i$ of a circle
attached to the vertex $i$.
Given $(M, \mathcal{T}, I)$, we assign each edge $\{ij\}$ the length
\begin{equation}\label{Euclidean length introduction}
l_{ij}=\sqrt{r_i^2+r_j^2+2r_ir_jI_{ij}}
\end{equation}
for Euclidean background geometry and
\begin{equation}\label{definition of length of edge for hyperbolic background 2}
l_{ij}=\cosh^{-1}(\cosh (r_i)\cosh(r_j)+I_{ij}\sinh(r_i)\sinh(r_j))
\end{equation}
for hyperbolic background geometry,
where $I_{ij}$ is the Euclidean and hyperbolic inversive distance of
the two circles centered at $v_i$ and $v_j$ with radii $r_i$ and $r_j$ respectively.
Note that the length $l_{ij}$ in (\ref{Euclidean length introduction}) and (\ref{definition of length of edge for hyperbolic background 2}) is well-defined for all $r_i>0, r_j>0$ under the condition $I_{ij}>-1$.
If $I_{ij}\in(-1,0)$, the two circles attached to the vertices $i$ and $j$ intersect with an obtuse angle.
If $I_{ij}\in[0,1]$, the two circles intersect with a non-obtuse angle. We can take $I_{ij}=\cos \Phi_{ij}$
with $\Phi_{ij}\in [0, \frac{\pi}{2}]$ and then the inversive distance circle packing is reduced to
Thurston's circle packing.
If $I_{ij}\in(1,+\infty)$, the two circles attached to the vertices $i$ and $j$ are disjoint.
See Figure \ref{Inversive_distance_circle_packing} for possible arrangements of the circles.
Guo \cite{Guo} and Luo \cite{L3} systematically studied the rigidity of inversive distance circle packing metrics
for nonnegative inversive distance $I\geq 0$, i.e. $I_{ij}\geq 0$ for every edge $\{ij\}\in E$.
In this paper, we focus on the case that $I>-1$.

The following is our main result, which solves Bowers and Stephenson's rigidity conjecture for inversive distance in $(-1, +\infty)$.

\begin{theorem}\label{main theorem global rigidity for Curvature}
Given a closed triangulated surface $(M, \mathcal{T}, I)$ with inversive distance $I: E\rightarrow (-1, +\infty)$ satisfying
\begin{equation}\label{nonnegative condition}
I_{ij}+I_{ik}I_{jk}\geq 0, I_{ik}+I_{ij}I_{jk}\geq 0, I_{jk}+I_{ij}I_{ik}\geq 0
\end{equation}
for any topological triangle $\triangle ijk\in F$.
\begin{description}
  \item[(1)] A Euclidean inversive distance circle packing on $(M, \mathcal{T})$ is determined by its combinatorial curvature $K: V\rightarrow \mathbb{R}$ up to scaling.
  \item[(2)] A hyperbolic inversive distance circle packing on $(M, \mathcal{T})$ is determined by its combinatorial curvature $K: V\rightarrow \mathbb{R}$.
\end{description}
\end{theorem}

\begin{remark}
For $I\in [0,1]$, the above result was Andreev and Thurston's rigidity for
circle packing with intersection angles in $[0,\frac{\pi}{2}]$. For $I\in (-1,1]$, the above result was the rigidity for
circle packing with intersection angles in $[0,\pi)$ recently obtained by Zhou \cite{Z}.
For $I\geq 0$, the above result was the rigidity for
inversive distance circle packing obtained by Guo \cite{Guo} and Luo \cite{L3}.
Our result unifies these results and allows the inversive distances to take values in a larger domain.
\end{remark}

\begin{remark}
It is interesting to note that in Theorem \ref{main theorem global rigidity for Curvature}, for a topological triangle $\triangle ijk\in F$,
if one of $I_{ij}, I_{ik}, I_{jk}$ is negative,
the other two must be nonnegative.
So at most one of $I_{ij}, I_{ik}, I_{jk}$ is negative.
\end{remark}

We further extend the rigidity to combinatorial $\alpha$-curvature introduced in \cite{GJ4, GX2, GX4, GX3, GX5, GX6}, which is defined as
$$R_{\alpha, i}=\frac{K_i}{s_i^\alpha}$$
for $\alpha\in \mathbb{R}$,
where $s_i=r_i$ for the Euclidean background geometry and $s_i=\tanh \frac{r_i}{2}$ for the hyperbolic background geometry.

\begin{theorem}\label{main theorem global rigidity for alpha curvature}
Given a closed triangulated surface $(M, \mathcal{T}, I)$ with inversive distance $I: E\rightarrow (-1, +\infty)$ satisfying
$$I_{ij}+I_{ik}I_{jk}\geq 0, I_{ik}+I_{ij}I_{jk}\geq 0, I_{jk}+I_{ij}I_{ik}\geq 0$$
for any topological triangle $\triangle ijk\in F$. $\overline{R}$
is a given function defined on the vertices  of  $(M, \mathcal{T})$.
\begin{description}
  \item[(1)]  If $\alpha\overline{R}\equiv0$, there exists at most one Euclidean inversive distance circle
  packing metric
  with combinatorial $\alpha$-curvature $\overline{R}$ up to scaling.
  If $\alpha\overline{R}\leq0$ and $\alpha\overline{R}\not\equiv0$,
  there exists at most one Euclidean inversive distance circle packing
  metric with combinatorial $\alpha$-curvature $\overline{R}$.
  \item[(2)] If $\alpha\overline{R}\leq 0$,
  there exists at most one hyperbolic inversive distance packing metric
  with combinatorial $\alpha$-curvature $\overline{R}$.
\end{description}
\end{theorem}

\subsection{Plan of paper}
The paper is organized as follows.
In Section \ref{section 3}, we study the Euclidean inversive distance circle packing metrics and prove Theorem
\ref{main theorem global rigidity for Curvature} and \ref{main theorem global rigidity for alpha curvature}
for the Euclidean background geometry.
In Section \ref{section 4}, we study the hyperbolic inversive distance circle packing metrics and prove Theorem
\ref{main theorem global rigidity for Curvature}, \ref{main theorem global rigidity for alpha curvature}
for the hyperbolic background geometry.

\section{Euclidean inversive distance circle packings}\label{section 3}

\subsection{Admissible space of Euclidean inversive distance circle packing metrics for a single triangle}

Given a weighted triangulated surface $(M, \mathcal{T}, I)$ with weight $I>-1$. Suppose
$\triangle ijk$ is a topological triangle in $F$. Here and in the following, to simplify notations,
when we are discussing a triangle $\triangle ijk$, we use $l_i$ to denote the length
of the edge $\{jk\}$ and use $I_i$ to denote the inversive distance of the two circles at the vertices $j$ and $k$.
In the Euclidean background geometry, the length $l_{i}$ of the edge $\{jk\}$ is then defined by
\begin{equation}\label{definition of Euclidean length}
l_{i}=\sqrt{r_j^2+r_k^2+2r_jr_kI_{i}}.
\end{equation}
For $I>-1$,  in order that the lengths $l_{i}, l_{j}, l_{k}$ for
$\Delta ijk\in F$ satisfy the triangle inequalities, there are some restrictions on the radii.
Denote the admissible space of the radius vectors for a face $\Delta ijk\in F$ as
\begin{equation}\label{admissible space}
\Omega^E_{ijk}:=\{(r_i, r_j, r_k)\in \mathbb{R}^3_{>0}|l_{i}<l_{j}+l_{k}, l_{j}<l_{i}+l_{k}, l_{k}<l_{i}+l_{j} \}.
\end{equation}
In the case of  $I\in [0,1]$, as noted by Thurston \cite{T1}, $\Omega^E_{ijk}=\mathbb{R}_{>0}^3$.
However, in general, $\Omega^E_{ijk}\neq\mathbb{R}_{>0}^3$ for $I\in(-1,+\infty)$.
It is proved \cite{Guo} that the admissible space $\Omega^E_{ijk}$ for $I\geq 0$ is a simply connected open subset of $\mathbb{R}_{>0}^3$
and $\Omega^E_{ijk}$ may not be convex.
Set
\begin{equation}
\Omega^E=\cap_{\Delta ijk\in F}\Omega^E_{ijk}
\end{equation}
to be the space of admissible radius function on the surface. $\Omega^E$ is obviously an open subset of $\mathbb{R}^N_{>0}$.
Every $r\in \Omega$ is called an inversive distance circle packing metric.

As noted in \cite{Guo}, in order that the edge lengths $l_{i}, l_{j}, l_{k}$ satisfy the triangle inequalities,
we just need
\begin{equation}\label{expansion of equiv tri inequ}
\begin{aligned}
0<&(l_{i}+l_{j}+l_{k})(l_{i}+l_{j}-l_{k})(l_{i}+l_{k}-l_{j})(l_{j}+l_{k}-l_{i})\\
=&4l_{i}^2l_{k}^2-(l_{i}^2+l_{k}^2-l_{j}^2)^2\\
=&2l_{i}^2l_{j}^2+2l_{i}^2l_{k}^2+2l_{j}^2l_{k}^2-l_{i}^4-l_{j}^4-l_{k}^4.
\end{aligned}
\end{equation}
Substituting the definition of edge length (\ref{definition of Euclidean length}) in the Euclidean background geometry
into (\ref{expansion of equiv tri inequ}),
by direct calculations, we have
\begin{equation*}
\begin{aligned}
 &\frac{1}{4}(l_{i}+l_{j}+l_{k})(l_{i}+l_{j}-l_{k})(l_{i}+l_{k}-l_{j})(l_{j}+l_{k}-l_{i})\\
=&r_i^2r_j^2(1-I^2_{k})+r_i^2r_k^2(1-I^2_{j})+r_j^2r_k^2(1-I^2_{i})\\
 &+2r_i^2r_jr_k(I_{i}+I_{j}I_{k})+2r_ir_j^2r_k(I_{j}+I_{i}I_{k})
 +2r_ir_jr_k^2(I_{k}+I_{i}I_{j})>0.
\end{aligned}
\end{equation*}
Denote
\begin{equation}\label{notation of gamma_ijk}
\gamma_{ijk}:=I_i+I_jI_k, \gamma_{jik}:=I_j+I_iI_k, \gamma_{kij}:=I_k+I_iI_j,
\end{equation}
then we have the following result on Euclidean triangle inequalities.
\begin{lemma} [\cite{Guo}]\label{Euclidean triangle inequality}
Suppose $(M, \mathcal{T}, I)$ is a weighted triangulated surface with weight $I>-1$ and
$\triangle ijk$ is a topological triangle in $F$.
The edge lengths $l_i, l_j, l_k$ defined by (\ref{definition of Euclidean length}) satisfy the triangle inequalities if and only if
\begin{equation}\label{equivalent triangle ineuqality}
\begin{aligned}
r_i^2r_j^2(1-I^2_{k})+r_i^2r_k^2(1-I^2_{j})+r_j^2r_k^2(1-I^2_{i})+2r_i^2r_jr_k\gamma_{ijk}+2r_ir_j^2r_k\gamma_{jik}
 +2r_ir_jr_k^2\gamma_{kij}>0.
\end{aligned}
\end{equation}
\end{lemma}

We have the following direct corollary obtained in \cite{Z} by Lemma  \ref{Euclidean triangle inequality}.
\begin{corollary}\label{Zhou's triangle condition}
If $I_i, I_j, I_k\in (-1, 1]$ and $\gamma_{ijk}\geq 0$, $\gamma_{jik}\geq 0$, $\gamma_{kij}\geq 0$,
then the triangle inequalities are satisfied for any $(r_i, r_j, r_k)\in \mathbb{R}^3_{>0}$.
\end{corollary}

\begin{remark}
Specially,  if $I_i=\cos \Phi_{i}, I_j=\cos \Phi_{j}, I_k=\cos \Phi_{k}$
with $\Phi_{i}, \Phi_{j}, \Phi_{k}\in [0, \frac{\pi}{2}]$,
then we have $I_i, I_j, I_k\in (-1, 1]$ and $\gamma_{ijk}\geq 0$, $\gamma_{jik}\geq 0$, $\gamma_{kij}\geq 0$.
So the triangle inequalities are satisfied for all radius vectors in $\mathbb{R}^3_{>0}$,
which was obtained by Thurston in \cite{T1}.
However, if we only require $\Phi_{i}, \Phi_{j}, \Phi_{k}\in [0, \pi)$, then (\ref{equivalent triangle ineuqality}) is equivalent to
\begin{equation*}
\begin{aligned}
&r_i^2r_j^2\sin^2\Phi_k+r_i^2r_k^2\sin^2\Phi_j+r_j^2r_k^2\sin^2\Phi_i+2r_i^2r_jr_k(\cos \Phi_i+\cos\Phi_j\cos\Phi_k)\\
&+2r_ir_j^2r_k(\cos \Phi_j+\cos\Phi_i\cos\Phi_k)
 +2r_ir_jr_k^2(\cos \Phi_k+\cos\Phi_i\cos\Phi_j)>0.
\end{aligned}
\end{equation*}
Specially, if $\Phi_i+\Phi_j\leq \pi, \Phi_i+\Phi_k\leq \pi, \Phi_j+\Phi_k\leq \pi$ \cite{Z}, or
$\Phi_i=\Phi_j\in [0, \frac{\pi}{2}]$ \cite{Z}, or $\Phi_i=\Phi_j=\Phi_k\in [0,\pi)$,
the triangle inequalities are satisfied.
\end{remark}

By Lemma \ref{Euclidean triangle inequality}, the admissible space
$\Omega^E_{ijk}$
for the topological triangle $\triangle ijk\in F$
may not be the whole space $\mathbb{R}^3_{>0}$. Furthermore, it is not always convex for all $I_i, I_j, I_k\in (-1, +\infty)$.
However, we have the following useful lemma on the structure of the admissible space $\Omega^E_{ijk}$.

\begin{lemma}\label{Euclidean admissible space simply connected}
Given a weighted triangulated surface $(M, \mathcal{T}, I)$ with $I>-1$.
For a topological triangle $\triangle ijk\in F$, if
\begin{equation}\label{gamma inequality condition}
\gamma_{ijk}\geq 0, \gamma_{jik}\geq 0, \gamma_{kij}\geq 0,
\end{equation}
then the admissible space $\Omega^E_{ijk}$ is a simply connected open subset of $\mathbb{R}^3_{>0}$.
Furthermore, for each connected component $V$ of $\mathbb{R}^3_{>0}\setminus \Omega^E_{ijk}$, the
intersection $V\cap \overline{\Omega}^E_{ijk}$ is a connected component of $\overline{\Omega}^E_{ijk}\setminus \Omega^E_{ijk}$,
on which $\theta_i$ is a constant function.
\end{lemma}

\proof
Define
\begin{equation*}
\begin{aligned}
F: \mathbb{R}^3_{>0}&\rightarrow \mathbb{R}^3_{>0}\\
(r_i, r_j, r_k)&\mapsto (r_j^2+r_k^2+2r_jr_kI_i, r_i^2+r_k^2+2r_ir_kI_j,r_i^2+r_j^2+2r_ir_jI_k)
\end{aligned}
\end{equation*}
and
\begin{equation*}
\begin{aligned}
G:\mathbb{R}^3_{>0}&\rightarrow \mathbb{R}^3_{>0}\\
(l_i, l_j, l_k)&\mapsto (l_i^2, l_j^2, l_k^2),
\end{aligned}
\end{equation*}
then $G$ is a diffeomorphism of $\mathbb{R}^3_{>0}$ and
$H=G^{-1}\circ F$ is the map sending $(r_i, r_j, r_k)$ to $(l_i, l_j, l_k)$.

We first prove that $H$ is injective. To prove this, we just need to prove that $F$ is injective.
Note that
\begin{equation*}
\begin{aligned}
\frac{\partial (F_i, F_j, F_k)}{\partial (r_i, r_j, r_k)}
=2\left(
   \begin{array}{ccc}
     0 & r_j+r_kI_i & r_k+r_jI_i \\
     r_i+r_kI_j & 0 & r_k+r_iI_j \\
     r_i+r_jI_k & r_j+r_iI_k & 0 \\
   \end{array}
 \right),
\end{aligned}
\end{equation*}
which implies that
\begin{equation*}
\begin{aligned}
\left|\frac{\partial (F_i, F_j, F_k)}{\partial (r_i, r_j, r_k)}\right|
=&8(r_j+r_kI_i)(r_k+r_iI_j)(r_k+r_iI_j)+8(r_k+r_jI_i)(r_i+r_kI_j)(r_j+r_iI_k)\\
=&16r_ir_jr_k(1+I_iI_jI_k)+8r_i\gamma_{ijk}(r_j^2+r_k^2)+8r_j\gamma_{jik}(r_i^2+r_k^2)+8r_k\gamma_{kij}(r_i^2+r_j^2).
\end{aligned}
\end{equation*}
By the condition (\ref{gamma inequality condition}) and the Cauchy inequality, we have
\begin{equation*}
\begin{aligned}
\left|\frac{\partial (F_i, F_j, F_k)}{\partial (r_i, r_j, r_k)}\right|
\geq& 16 r_ir_jr_k(1+I_iI_jI_k+\gamma_{ijk}+\gamma_{jik}+\gamma_{kij})\\
=&16 r_ir_jr_k(1+I_i)(1+I_j)(1+I_k).
\end{aligned}
\end{equation*}
By the condition that $I_i, I_j, I_k\in (-1,+\infty)$, we have $\left|\frac{\partial (F_i, F_j, F_k)}{\partial (r_i, r_j, r_k)}\right|>0$
for any $r\in \mathbb{R}^3_{>0}$.
If there are $r=(r_i, r_j, r_k)\in \mathbb{R}^3_{>0}$ and $r'=(r'_i, r'_j, r'_k)\in \mathbb{R}^3_{>0}$
satisfying $F(r)=F(r')$, then we have
\begin{equation*}
\begin{aligned}
0=F(r)-F(r')=\frac{\partial (F_i, F_j, F_k)}{\partial (r_i, r_j, r_k)}\big|_{r+\theta(r-r')}\cdot(r-r')^T,\ \ \theta\in(0,1),
\end{aligned}
\end{equation*}
which implies $r=r'$ by the nondegeneracy of $\frac{\partial (F_i, F_j, F_k)}{\partial (r_i, r_j, r_k)}$ on $\mathbb{R}^3_{>0}$.
So the map $F$ is injective on $\mathbb{R}^3_{>0}$, which implies that $H$ is injective on $\mathbb{R}^3_{>0}$.

Note that
\begin{equation*}
\begin{aligned}
F_{i}=&r_j^2+r_k^2+2r_jr_kI_i\geq 2r_jr_k(1+I_i),\\
F_{j}=&r_i^2+r_k^2+2r_ir_jI_k\geq 2r_ir_k(1+I_j),\\
F_{k}=&r_i^2+r_j^2+2r_ir_jI_k\geq 2r_ir_j(1+I_k).
\end{aligned}
\end{equation*}
By the condition that $I_i, I_j, I_k\in (-1, +\infty)$, if $F$ is bounded, we have $r_ir_j$, $r_ir_k$, $r_jr_k$ are bounded, which implies that
$r_i^2+r_j^2, r_i^2+r_k^2, r_j^2+r_k^2$ are bounded.
Similarly, we have $F_i\leq (1+|I_i|)(r_j^2+r_j^2)$.
This implies that $F$ is a proper map from
$\mathbb{R}^3_{>0}$ to $\mathbb{R}^3_{>0}$. By the invariance of domain theorem, we have
$F$ is a diffeomorphism between $\mathbb{R}^3_{>0}$ and $F(\mathbb{R}^3_{>0})$.
And then $H$ is a diffeomorphism between $\mathbb{R}^3_{>0}$ and $H(\mathbb{R}^3_{>0})$.

Set
$$\mathcal{L}=\{(l_i,l_j,l_k)|l_i+l_j>l_k, l_i+l_k>l_j, l_j+l_k>l_i\},$$
then $\Omega^E_{ijk}=H^{-1}(H(\mathbb{R}^3_{>0})\cap \mathcal{L})$. To prove that $\Omega^E_{ijk}$
is simply connected, we just need to prove that $H(\mathbb{R}^3_{>0})\cap \mathcal{L}$ is a cone.
Note that $\mathcal{L}$ is a cone in $\mathbb{R}^3_{>0}$ bounded by three planes
\begin{equation*}
\begin{aligned}
L_i=&\{(l_{i}, l_j, l_k)\in \mathbb{R}^3_{>0}|l_i=l_j+l_k\},\\
L_j=&\{(l_{i}, l_j, l_k)\in \mathbb{R}^3_{>0}|l_j=l_i+l_k\},\\
L_k=&\{(l_{i}, l_j, l_k)\in \mathbb{R}^3_{>0}|l_k=l_i+l_j\}.\\
\end{aligned}
\end{equation*}
Note that $H$ is a diffeomorphism between $\mathbb{R}^3_{>0}$ and $H(\mathbb{R}^3_{>0})$,
$H(\mathbb{R}^3_{>0})$ is a cone bounded by three quadratic surfaces
\begin{equation*}
\begin{aligned}
\Sigma_i=&\{(l_{i}, l_j, l_k)\in \mathbb{R}^3_{>0}|l_i^2=l_j^2+l_k^2+2l_jl_kI_i\},\\
\Sigma_i=&\{(l_{i}, l_j, l_k)\in \mathbb{R}^3_{>0}|l_j^2=l_i^2+l_k^2+2l_il_kI_j\},\\
\Sigma_i=&\{(l_{i}, l_j, l_k)\in \mathbb{R}^3_{>0}|l_k^2=l_i^2+l_j^2+2l_il_jI_k\}.
\end{aligned}
\end{equation*}
In fact, if $r_i=0$, then $l_j=r_k, l_k=r_j$ and
$l_i^2=r_j^2+r_k^2+2r_jr_kI_i=l_j^2+l_k^2+2l_jl_kI_i.$
$\Sigma_i$ is in fact the image of $r_i=0$ under $H$.
By the diffeomorphism of $H$, $\Sigma_i$, $\Sigma_j$, $\Sigma_k$ are mutually disjoint.
Furthermore, if $I_i\in (-1, 1]$, we have
$(l_j-l_k)^2<l_i^2\leq(l_j+l_k)^2$ on $\Sigma_i$. And if $I_i\in (1, +\infty)$, we have $l_i^2>(l_j+l_k)^2$ on $\Sigma_i$.
This implies that $\Sigma_i\subset \overline{\mathcal{L}}$ if $I_i\in (-1, 1]$ and $\Sigma_i\cap \mathcal{L}=\emptyset$ if $I_i\in (1, +\infty)$.
Similar results hold for $\Sigma_j$ and $\Sigma_k$.

To prove that $H(\mathbb{R}^3_{>0})\cap \mathcal{L}$ is a cone, we just need to consider the following cases by the symmetry between $i,j,k$.

If $I_i, I_j, I_k\in (-1, 1]$, $H(\mathbb{R}^3_{>0})\cap \mathcal{L}$ is a cone bounded by $\Sigma_i, \Sigma_j, \Sigma_k$ and
$H(\mathbb{R}^3_{>0})\cap \mathcal{L}=H(\mathbb{R}^3_{>0})$.

If $I_i, I_j\in (-1, 1]$ and $I_k\in(1,+\infty)$, $H(\mathbb{R}^3_{>0})\cap \mathcal{L}$ is a cone bounded by $\Sigma_i, \Sigma_j$ and $L_k$.

If $I_i\in (-1, 1]$ and $I_j, I_k\in(1,+\infty)$, $H(\mathbb{R}^3_{>0})\cap \mathcal{L}$ is a cone bounded by $\Sigma_i$, $L_j$ and $L_k$.

If $I_i, I_j, I_k\in(1,+\infty)$, $H(\mathbb{R}^3_{>0})\cap \mathcal{L}$ is a cone bounded by $L_i$, $L_j$ and $L_k$.
In this case, $H(\mathbb{R}^3_{>0})\cap \mathcal{L}=\mathcal{L}$.

For any case, $H(\mathbb{R}^3_{>0})\cap \mathcal{L}$ is a cone in $\mathbb{R}^3_{>0}$.
By the fact that $H$ is a diffeomorphism between $\mathbb{R}^3_{>0}$ and $H(\mathbb{R}^3_{>0})$,
we have the admissible space $\Omega^E_{ijk}=H^{-1}(H(\mathbb{R}^3_{>0})\cap \mathcal{L})$ is simply connected.

By the analysis above,
if $H(\mathbb{R}^3_{>0})\subset \mathcal{L}$, then $\Omega^E_{ijk}=H^{-1}(H(\mathbb{R}^3_{>0})\cap\mathcal{L})=\mathbb{R}^3_{>0}$.
If $H(\mathbb{R}^3_{>0})\setminus\mathcal{L}\neq \emptyset$, then $\Omega^E_{ijk}$ is a proper subset of $\mathbb{R}^3_{>0}$.
If $I_i>1$, the boundary component $\Sigma_i=\{l_i^2=l_j^2+l_k^2+2l_jl_kI_i\}$ is out of the set $\mathcal{L}$.
By the fact that $\Omega^E_{ijk}=H^{-1}(H(\mathbb{R}^3_{>0})\cap\mathcal{L})$ and
$H:\mathbb{R}^3_{>0}\rightarrow H(\mathbb{R}^3_{>0})$ is a diffeomorphism,
we have $H^{-1}(L_i)$ is a connected boundary component of $\Omega^E_{ijk}$, on which $\theta_i=\pi, \theta_j=\theta_k=0$.
This completes the proof of the lemma.
\qed

\begin{corollary}\label{Euclidean triangle extension lemma}
For a topological triangle $\triangle ijk\in F$ with inversive distance $I>-1$ and
$
\gamma_{ijk}\geq 0, \gamma_{jik}\geq 0, \gamma_{kij}\geq 0,
$
the functions $\theta_i, \theta_j, \theta_k$ defined on $\Omega^E_{ijk}$
could be continuously extended by constant to
$\widetilde{\theta}_i, \widetilde{\theta}_j, \widetilde{\theta}_k$  defined on $\mathbb{R}^3_{>0}$.
\end{corollary}

\begin{remark}\label{remark Euclidean triangle inequ}
\begin{description}
  \item[(1)] If $I_i, I_j, I_k\in [0, +\infty)$, obviously we have $\gamma_{ijk}\geq 0$, $\gamma_{jik}\geq 0$, $\gamma_{kij}\geq 0$. So
  Lemma \ref{Euclidean admissible space simply connected} generalizes Lemma 3 in \cite{Guo} obtained by Guo.
  \item[(2)] If $I_i, I_j, I_k\in (-1, 1]$ and $\gamma_{ijk}\geq 0$, $\gamma_{jik}\geq 0$, $\gamma_{kij}\geq 0$,
  by the proof of Lemma \ref{Euclidean admissible space simply connected}, $\Omega^E_{ijk}=\mathbb{R}^3_{>0}$, which is obtained by Zhou \cite{Z} .
  \item[(3)] The condition $I_i, I_j, I_k\in (-1, +\infty)$ and $\gamma_{ijk}\geq 0$, $\gamma_{jik}\geq 0$, $\gamma_{kij}\geq 0$ contains more
  cases, for example, $I_i=-\frac{1}{2}$, $I_j=1$ and $I_k=2$, in which case the admissible space $\Omega^E_{ijk}$ is still simply connected.
\end{description}
\end{remark}

\subsection{Infinitesimal rigidity of Euclidean inversive distance circle packings}

Set $u_i=\ln r_i$, then we have $\mathcal{U}^E_{ijk}:=\ln (\Omega^E_{ijk})$ is a simply connected subset of $\mathbb{R}^3$
by Lemma \ref{Euclidean admissible space simply connected}.
If $(r_i, r_j, r_k)\in \Omega^E_{ijk}$, $l_i, l_j, l_k$ satisfy the triangle inequalities and forms a Euclidean triangle.
Denote the inner angle at the vertex $i$ as $\theta_i$. We have the following useful lemma.

\begin{lemma}\label{Euclidean symmetry lemma}
For any topological triangle $\triangle ijk\in F$, we have
\begin{equation}\label{derivative of theta}
\frac{\partial \theta_i}{\partial u_j}=\frac{\partial \theta_j}{\partial u_i}
=\frac{1}{Al_{k}^2}\left[r_i^2r_j^2(1-I_k^2)+r_i^2r_jr_k\gamma_{ijk}+r_ir_j^2r_k\gamma_{jik}\right]
\end{equation}
on $\mathcal{U}^E_{ijk}$, where $A=l_jl_k\sin \theta_i$.
\end{lemma}
\proof
By the cosine law, we have
$l_i^2=l_j^2+l_k^2-2l_jl_k\cos\theta_i.$
Taking the derivative with respect to $l_i$, we have
$\frac{\partial \theta_i}{\partial l_i}=\frac{l_i}{A},$
where $A=l_jl_k\sin \theta_i$ is two times of the area of $\triangle ijk$.
Similarly, we have
$\frac{\partial \theta_i}{\partial l_j}=\frac{-l_i\cos\theta_k}{A},
\frac{\partial \theta_i}{\partial l_k}=\frac{-l_i\cos\theta_j}{A}.$
By the definition of $l_i, l_j, l_k$, we have
$$\frac{\partial l_i}{\partial r_j}=\frac{r_j+r_kI_i}{l_i},
\frac{\partial l_j}{\partial r_j}=0,
\frac{\partial l_k}{\partial r_j}=\frac{r_j+r_iI_k}{l_k}.$$
Then
\begin{equation*}
\begin{aligned}
\frac{\partial \theta_i}{\partial u_j}
=&r_j\frac{\partial \theta_i}{\partial r_j}\\
=&r_j(\frac{\partial \theta_i}{\partial l_i}\frac{\partial l_i}{\partial r_j}+\frac{\partial \theta_i}{\partial l_k}\frac{\partial l_k}{\partial r_j}
)\\
=&r_j\left[\frac{r_j+r_kI_i}{A}-\frac{l_i\cos\theta_j(r_j+r_iI_k)}{Al_k}\right]\\
=&\frac{1}{Al_k}\left[l_k(r_j^2-r_jr_kI_i)-\frac{l_i^2+l_k^2-l_j^2}{2l_k}(r_j^2+r_ir_jI_k)\right]\\
=&\frac{1}{Al_{k}^2}\left[r_i^2r_j^2(1-I_k^2)+r_i^2r_jr_k\gamma_{ijk}+r_ir_j^2r_k\gamma_{jik}\right],
\end{aligned}
\end{equation*}
where the cosine law is used in the third line and the definition of the length (\ref{definition of Euclidean length})
is used in the fourth line.
This also implies $\frac{\partial \theta_i}{\partial u_j}=\frac{\partial \theta_j}{\partial u_i}$. \qed

\begin{remark}\label{Euclidean mononicity of angle}
The equation $\frac{\partial \theta_i}{\partial u_j}=\frac{\partial \theta_j}{\partial u_i}$ has been
obtained under different conditions in \cite{CL1,DV,Guo} and
the formulas for $\frac{\partial \theta_i}{\partial l_j}$ and $\frac{\partial \theta_i}{\partial l_i}$ was obtained
by Chow and Luo \cite{CL1}.
In general, for $I_{i}, I_j, I_k\in (-1, +\infty)$, $\frac{\partial \theta_i}{\partial u_j}$ have no sign.
However, if $I_{i}, I_j, I_k\in (-1, 1]$ and $\gamma_{ijk}\geq 0$, $\gamma_{jik}\geq 0$, $\gamma_{kij}\geq 0$, by (\ref{derivative of theta}),
we have $\frac{\partial \theta_i}{\partial u_j}\geq 0$. Furthermore, $\frac{\partial \theta_i}{\partial u_j}=0$ if and only if
$I_k=1$ and $I_i+I_j=0$. Especially, if $I_i=\cos\Phi_{i}, I_j=\cos\Phi_{j}, I_{k}=\cos\Phi_{k}$
with $\Phi_{i},\Phi_{j},\Phi_{k}\in [0,\frac{\pi}{2}]$, we have
$\frac{\partial \theta_i}{\partial u_j}\geq 0$, and $\frac{\partial \theta_i}{\partial u_j}=0$ if and only if
$\Phi_{k}=0$ and $\Phi_{i}=\Phi_{j}=\frac{\pi}{2}$.
\end{remark}
\begin{remark}
Geometrically, the three circles at the vertices have a power center $O$.
It is known \cite{ZG,ZGZLYG} that $\frac{\partial \theta_i}{\partial u_j}=\frac{h_k}{l_k}$, where $h_k$ is the signed distance of
the power center $O$ to the edge $\{ij\}$, which is positive if $O$ is in the interior of the triangle $\triangle ijk$ and negative if
the power center $O$ is out of the triangle $\triangle ijk$. So under the condition
$I_{i}, I_j, I_k\in (-1, 1]$ and $\gamma_{ijk}\geq 0$, $\gamma_{jik}\geq 0$, $\gamma_{kij}\geq 0$,
the power center $O$ is in the triangle $\triangle ijk$.
\end{remark}

Lemma \ref{Euclidean symmetry lemma} shows that the matrix
\begin{equation*}
\begin{aligned}
\Lambda^E_{ijk}:=\frac{\partial (\theta_i, \theta_j, \theta_k)}{\partial (u_i, u_j, u_k)}
=\left(
   \begin{array}{ccc}
     \frac{\partial \theta_i}{\partial u_i} & \frac{\partial \theta_i}{\partial u_j} & \frac{\partial \theta_i}{\partial u_k} \\
     \frac{\partial \theta_j}{\partial u_i} & \frac{\partial \theta_j}{\partial u_j} & \frac{\partial \theta_j}{\partial u_k} \\
     \frac{\partial \theta_k}{\partial u_i} & \frac{\partial \theta_k}{\partial u_j} & \frac{\partial \theta_k}{\partial u_k} \\
   \end{array}
 \right)
\end{aligned}
\end{equation*}
is symmetric on $\mathcal{U}^E_{ijk}$. For the matrix $\Lambda^E_{ijk}$, we have the following useful property.

\begin{lemma}\label{Euclidean negative semidefinite}
For any topological triangle $\triangle ijk\in F$ with inversive distance $I_{i}, I_j, I_k\in (-1, +\infty)$ and
$\gamma_{ijk}\geq 0$, $\gamma_{jik}\geq 0$, $\gamma_{kij}\geq 0$,
the matrix $\Lambda^E_{ijk}$ is negative semi-definite with rank 2 and kernel $\{t(1,1,1)^T|t\in \mathbb{R}\}$ on $\mathcal{U}^E_{ijk}$.
\end{lemma}
\proof
The proof is parallel to that of Lemma 6 in \cite{Guo} with some modifications.
By the calculations in Lemma \ref{Euclidean symmetry lemma}, for a triangle $\triangle ijk\in F$, we have
\begin{equation*}
\begin{aligned}
\left(
  \begin{array}{c}
    d\theta_i \\
    d\theta_j \\
    d\theta_k \\
  \end{array}
\right)
=&-\frac{1}{A} \left(
                \begin{array}{ccc}
                  l_i & 0 & 0 \\
                  0 & l_j & 0 \\
                  0 & 0 & l_k \\
                \end{array}
              \right)
              \left(
                \begin{array}{ccc}
                  -1 & \cos\theta_k & \cos\theta_j \\
                  \cos\theta_k & -1 & \cos\theta_i \\
                  \cos\theta_j & \cos\theta_i & -1 \\
                \end{array}
              \right)\\
  &\times\left(
     \begin{array}{ccc}
       0 & \frac{l_i^2+r_j^2-r_k^2}{2l_ir_j} & \frac{l_i^2+r_k^2-r_j^2}{2l_ir_k} \\
       \frac{l_j^2+r_i^2-r_k^2}{2l_jr_i} & 0 & \frac{l_j^2+r_k^2-r_i^2}{2l_jr_k} \\
       \frac{l_k^2+r_i^2-r_j^2}{2l_kr_i} & \frac{l_k^2+r_j^2-r_i^2}{2l_kr_i} & 0 \\
     \end{array}
   \right)
   \left(
     \begin{array}{ccc}
       r_i & 0 & 0 \\
       0 & r_j & 0 \\
       0 & 0 & r_k \\
     \end{array}
   \right)
   \left(
     \begin{array}{c}
       du_i \\
       du_j \\
       du_k \\
     \end{array}
   \right).
\end{aligned}
\end{equation*}
Write the above formula as
\begin{equation*}
\begin{aligned}
\left(
  \begin{array}{c}
    d\theta_i \\
    d\theta_j \\
    d\theta_k \\
  \end{array}
\right)
=-\frac{1}{A}N
   \left(
     \begin{array}{c}
       du_i \\
       du_j \\
       du_k \\
     \end{array}
   \right).
\end{aligned}
\end{equation*}
By the cosine law, we have
\begin{equation*}
\begin{aligned}
4N=&\left(
     \begin{array}{ccc}
       -2l_i^2 & l_i^2+l_j^2-l_k^2 & l_k^2+l_i^2-l_j^2 \\
       l_i^2+l_j^2-l_k^2 & -2l_j^2 & l_j^2+l_k^2-l_i^2 \\
       l_k^2+l_i^2-l_j^2 & l_j^2+l_k^2-l_i^2 & -2l_k^2 \\
     \end{array}
   \right)
      \left(
     \begin{array}{ccc}
       \frac{1}{l_i^2} & 0 & 0 \\
       0 & \frac{1}{l_j^2} & 0 \\
       0 & 0 & \frac{1}{l_k^2} \\
     \end{array}
   \right)\\
   &\times\left(
     \begin{array}{ccc}
       0 & l_i^2+r_j^2-r_k^2 & l_i^2+r_k^2-r_j^2 \\
       l_j^2+r_i^2-r_k^2 & 0 & l_j^2+r_k^2-r_i^2 \\
       l_k^2+r_i^2-r_j^2 & l_k^2+r_j^2-r_i^2 & 0 \\
     \end{array}
   \right)
\end{aligned}
\end{equation*}
By Lemma \ref{Euclidean symmetry lemma}, we have $4N$ is symmetric. Furthermore, note that
$\theta_i+\theta_j+\theta_k=\pi$, we have
$0=\frac{\partial \theta_i}{\partial u_i}+\frac{\partial \theta_j}{\partial u_i}+\frac{\partial \theta_k}{\partial u_i}
=\frac{\partial \theta_i}{\partial u_i}+\frac{\partial \theta_i}{\partial u_j}+\frac{\partial \theta_i}{\partial u_k}.$
Then we can write $4N$ as
\begin{equation*}
\begin{aligned}
4N=\left(
     \begin{array}{ccc}
       -A-B & A & B \\
       A & -A-C & C \\
       B & C & -B-C \\
     \end{array}
   \right).
\end{aligned}
\end{equation*}
To prove $\Lambda^E_{ijk}$ is negative semi-definite, we just need to prove that $4N$ is positive semi-definite.
By direct calculations, we have
\begin{equation*}
\begin{aligned}
|\lambda I-4N|
=&\left|
         \begin{array}{ccc}
           \lambda+A+B & -A & -B \\
           -A & \lambda+A+C & -C \\
           -B & -C & \lambda+B+C \\
         \end{array}
\right|\\
=&\lambda[\lambda^2+2(A+B+C)\lambda+3(AB+AC+BC)].
\end{aligned}
\end{equation*}
We want to prove that the equation
$$\lambda^2+2(A+B+C)\lambda+3(AB+AC+BC)=0$$
has two positive roots.
Note that for this quadratic equation, we have
$$\Delta=4(A+B+C)^2-12(AB+AC+BC)=4(A^2+B^2+C^2-AB-AC-BC)\geq 0,$$
so we just need to prove that $A+B+C<0$ and $AB+AC+BC>0$.

By direct calculations, we have
\begin{equation*}
\begin{aligned}
-2(A+B+C)=l_i^2+l_j^2+l_k^2+(l_j^2-l_k^2)\frac{r_j^2-r_k^2}{l_i^2}+(l_k^2-l_i^2)\frac{r_j^2-r_i^2}{l_j^2}+(l_i^2-l_j^2)\frac{r_i^2-r_j^2}{l_k^2}.
\end{aligned}
\end{equation*}
So $A+B+C<0$ is equivalent to
$$l_i^2+l_j^2+l_k^2+(l_j^2-l_k^2)\frac{r_j^2-r_k^2}{l_i^2}+(l_k^2-l_i^2)\frac{r_j^2-r_i^2}{l_j^2}+(l_i^2-l_j^2)\frac{r_i^2-r_j^2}{l_k^2}>0,$$
which is equivalent to
$$l_i^2l_j^2l_k^2(l_i^2+l_j^2+l_k^2)+l_i^2r_i^2(l_i^2l_j^2+l_i^2l_k^2-l_j^4-l_k^4)+l_j^2r_j^2(l_i^2l_j^2+l_j^2l_k^2-l_i^4-l_k^4)
+l_k^2r_k^2(l_i^2l_k^2+l_j^2l_k^2-l_i^4-l_j^4)>0.$$
Note that
\begin{equation*}
\begin{aligned}
&2[l_i^2l_j^2l_k^2(l_i^2+l_j^2+l_k^2)+l_i^2r_i^2(l_i^2l_j^2+l_i^2l_k^2-l_j^4-l_k^4)\\
&+l_j^2r_j^2(l_i^2l_j^2+l_j^2l_k^2-l_i^4-l_k^4)
+l_k^2r_k^2(l_i^2l_k^2+l_j^2l_k^2-l_i^4-l_j^4)]\\
=&2l_i^2l_j^2l_k^2(l_i^2+l_j^2+l_k^2)+l_i^2r_i^2(l_i^4-l_j^4-l_k^4-2l_j^2l_k^2)+l_j^2r_j^2(l_j^4-l_i^4-l_k^4-2l_i^2l_k^2)\\
&+l_k^2r_k^2(l_k^4-l_i^4-l_j^4-2l_i^2l_j^2)+(l_i^2r_i^2+l_j^2r_j^2+l_k^2r_k^2)(2l_i^2l_j^2+2l_i^2l_k^2+2l_j^2l_k^2-l_i^4-l_j^4-l_k^4).
\end{aligned}
\end{equation*}
By the triangle inequalities, we have
$$2l_i^2l_j^2+2l_i^2l_k^2+2l_j^2l_k^2-l_i^4-l_j^4-l_k^4>0$$
on $\Omega^E_{ijk}$.
So to prove $A+B+C<0$, we just need to prove
\begin{equation*}
\begin{aligned}
&2l_i^2l_j^2l_k^2(l_i^2+l_j^2+l_k^2)+l_i^2r_i^2(l_i^4-l_j^4-l_k^4-2l_j^2l_k^2)+l_j^2r_j^2(l_j^4-l_i^4-l_k^4-2l_i^2l_k^2)\\
&+l_k^2r_k^2(l_k^4-l_i^4-l_j^4-2l_i^2l_j^2)>0.
\end{aligned}
\end{equation*}
By direct calculations, we have
\begin{equation*}
\begin{aligned}
&2l_i^2l_j^2l_k^2(l_i^2+l_j^2+l_k^2)+l_i^2r_i^2(l_i^4-l_j^4-l_k^4-2l_j^2l_k^2)\\
&+l_j^2r_j^2(l_j^4-l_i^4-l_k^4-2l_i^2l_k^2)+l_k^2r_k^2(l_k^4-l_i^4-l_j^4-2l_i^2l_j^2)\\
=&4[r_i^2r_j^2r_k^2(1+I_i^2+I_j^2+I_k^2+4I_iI_jI_k)+r_i^2r_jr_k(I_i+I_jI_k)(r_j^2+r_k^2)\\
&+r_ir_j^2r_k(I_j+I_iI_k)(r_i^2+r_k^2)+r_ir_jr_k^2(I_k+I_iI_j)(r_i^2+r_j^2)]\\
\geq&4r_i^2r_j^2r_k^2(1+I_i^2+I_j^2+I_k^2+4I_iI_jI_k+2I_i+2I_jI_k+2I_j+2I_iI_k+2I_k+2I_iI_j)\\
=&4r_i^2r_j^2r_k^2[(1+I_i)(1+I_j)(1+I_k)+(1+I_i)\gamma_{ijk}+(1+I_j)\gamma_{jik}+(1+I_k)\gamma_{kij}]\\
>&0,
\end{aligned}
\end{equation*}
where the condition $I_i, I_j, I_k\in (-1, +\infty)$ and $\gamma_{ijk}=I_i+I_jI_k\geq 0, \gamma_{jik}=I_j+I_iI_k\geq 0, \gamma_{kij}=I_k+I_iI_j\geq 0$
is used.
So we have $A+B+C<0$.

For the term $AB+AC+BC$, by direct calculations, we have
\begin{equation*}
\begin{aligned}
&AB+AC+BC\\
=&\frac{1}{l_i^2l_j^2l_k^2}(2l_i^2l_j^2+2l_i^2l_k^2+2l_j^2l_k^2-l_i^4-l_j^4-l_k^4)\\
 &\times \left[(r_i^2-r_j^2)(r_k^2-r_i^2)l_i^2+(r_i^2-r_j^2)(r_j^2-r_k^2)l_j^2+(r_k^2-r_i^2)(r_j^2-r_k^2)l_k^2+l_i^2l_j^2l_k^2\right].
\end{aligned}
\end{equation*}
So by the triangle inequalities, $AB+AC+BC>0$ is equivalent to
$$(r_i^2-r_j^2)(r_k^2-r_i^2)l_i^2+(r_i^2-r_j^2)(r_j^2-r_k^2)l_j^2+(r_k^2-r_i^2)(r_j^2-r_k^2)l_k^2+l_i^2l_j^2l_k^2>0.$$
By direct calculations, combining with the condition $I_i, I_j, I_k\in (-1, +\infty)$ and $\gamma_{ijk}\geq 0, \gamma_{jik}\geq 0, \gamma_{kij}\geq 0$,
we have
\begin{equation*}
\begin{aligned}
&(r_i^2-r_j^2)(r_k^2-r_i^2)l_i^2+(r_i^2-r_j^2)(r_j^2-r_k^2)l_j^2+(r_k^2-r_i^2)(r_j^2-r_k^2)l_k^2+l_i^2l_j^2l_k^2\\
=&8r_i^2r_j^2r_k^2(1+I_iI_jI_k)+4r_i^2r_jr_k(I_i+I_jI_k)(r_j^2+r_k^2)\\
&+4r_ir_j^2r_k(I_j+I_iI_k)(r_i^2+r_k^2)+4r_ir_jr_k^2(I_k+I_iI_j)(r_i^2+r_j^2)\\
\geq& 8r_i^2r_j^2r_k^2 (1+I_iI_jI_k+I_i+I_jI_k+I_j+I_iI_k+I_k+I_iI_j)\\
=&8r_i^2r_j^2r_k^2 (1+I_i)(1+I_j)(1+I_k)\\
>&0.
\end{aligned}
\end{equation*}
So we have $AB+AC+BC>0$.
Then the matrix $\Lambda^E_{ijk}$ has a zero eigenvalue with eigenvector $(1,1,1)^T$ and two negative eigenvalues on $\mathcal{U}^E_{ijk}$.
\qed

Now suppose that for each topological face $\Delta ijk\in F$, the triangle inequalities are satisfied, i.e. $r\in \Omega^E$,
then the weighted triangulated surface
$(M, \mathcal{T}, I)$ could be taken as gluing many triangles along the edges coherently, which produces
a cone metric on the triangulated surface with singularities at $V$.
To describe the singularity at the vertex $i$, the classical discrete curvature is introduced, which is defined as
\begin{equation}\label{classical Gauss curv}
K_i=2\pi-\sum_{\triangle ijk \in F}\theta_i^{jk},
\end{equation}
where the sum is taken over all the triangles with $i$ as one of its vertices and
$\theta_i^{jk}$ is the inner angle of the triangle $\triangle ijk$ at the vertex $i$.
Lemma \ref{Euclidean negative semidefinite} has the following corollary.

\begin{corollary}\label{positivity of Lambda E}
Given a triangulated surface $(M, \mathcal{T})$ with inversive distance $I>-1$ and
$\gamma_{ijk}\geq 0, \gamma_{jik}\geq 0, \gamma_{kij}\geq 0$
for any topological triangle $\triangle ijk\in F$.
Then the matrix $\Lambda^E=\frac{\partial (K_1,\cdots, K_N)}{\partial(u_1,\cdots,u_N)}$
is symmetric and positive semi-definite with rank $N-1$ and kernel $\{t\textbf{1}|t\in \mathbb{R}\}$
on $\mathcal{U}^E$ for the Euclidean background geometry.
\end{corollary}
\proof
This follows from the fact that $\Lambda^E=-\sum_{\triangle ijk\in F}\Lambda^E_{ijk}$,
Lemma \ref{Euclidean symmetry lemma} and
Lemma \ref{Euclidean negative semidefinite}, where $\Lambda^E_{ijk}$ is extended by zeros to a $N\times N$
matrix so that $\Lambda^E_{ijk}$ acts on a vector $(v_1, \cdots, v_N)$ only on
the coordinates corresponding to vertices $v_i, v_j$ and $v_k$ in the triangle $\triangle ijk$.
\qed

\begin{remark}
Guo \cite{Guo} obtained a result paralleling to Corollary \ref{positivity of Lambda E} for nonnegative inversive distance.
\end{remark}

By Lemma \ref{Euclidean admissible space simply connected} and Lemma \ref{Euclidean symmetry lemma}, we can define an
energy function
$$\mathcal{E}_{ijk}(u)=\int_{u_0}^u\theta_idu_i+\theta_jdu_j+\theta_kdu_k$$
on $\mathcal{U}^E_{ijk}$. Lemma \ref{Euclidean negative semidefinite} ensures that
$\mathcal{E}_{ijk}$ is locally concave on $\mathcal{U}^E_{ijk}$.
Define the Ricci energy function as
\begin{equation}\label{local Euclid energy function}
\mathcal{E}(u)=-\sum_{\triangle ijk\in F}\mathcal{E}_{ijk}(u)+\int_{u_0}^u\sum_{i=1}^N(2\pi-\overline{K}_i)du_i,
\end{equation}
then $\nabla_u\mathcal{E}=K-\overline{K}$ and $\mathcal{E}(u)$ is locally convex on $\mathcal{U}^E=\cap_{\triangle ijk\in F} \mathcal{U}^E_{ijk}$.
The local convexity of $\mathcal{E}$ implies the infinitesimal rigidity of $K$ with respect to $u$, which is the infinitesimal rigidity of inversive
distance circle packings.

\subsection{Global rigidity of Euclidean inversive distance circle packings}

In this subsection, we shall prove the global rigidity of inversive distance circle packings under the condition $I>-1$ and
$\gamma_{ijk}\geq 0, \gamma_{jik}\geq 0, \gamma_{kij}\geq 0$
for any triangle $\triangle ijk\in F$.
We need to extend the energy function defined on $\mathcal{U}^E$ to be a convex function defined on $\mathbb{R}^3$.
Before going on, we recall the following definition and theorem of Luo in \cite{L3}.
\begin{definition}
A differential 1-form $w=\sum_{i=1}^n a_i(x)dx^i$ in an open set $U\subset \mathbb{R}^n$ is said to be continuous if
each $a_i(x)$ is continuous on $U$.  A differential 1-form $w$ is called closed if $\int_{\partial \tau}w=0$ for each
triangle $\tau\subset U$.
\end{definition}

\begin{theorem}[\cite{L3} Corollary 2.6]\label{Luo's convex extention}
Suppose $X\subset \mathds{R}^n$ is an open convex set and $A\subset X$ is an open subset of $X$ bounded by a $C^1$
smooth codimension-1 submanifold in $X$. If $w=\sum_{i=1}^na_i(x)dx_i$ is a continuous closed 1-form on $A$ so that
$F(x)=\int_a^x w$ is locally convex on $A$ and each $a_i$ can be extended continuous to $X$ by constant functions to a
function $\widetilde{a}_i$ on $X$, then  $\widetilde{F}(x)=\int_a^x\sum_{i=1}^n\widetilde{a}_i(x)dx_i$ is a $C^1$-smooth
convex function on $X$ extending $F$.
\end{theorem}

Combining Lemma \ref{Euclidean admissible space simply connected}, Corollary \ref{Euclidean triangle extension lemma}
and Theorem \ref{Luo's convex extention}, we have the following useful lemma.

\begin{lemma}\label{Euclidean extension lemma}
For any triangle $\triangle ijk\in F$ with inversive distance $I>-1$ and
$$\gamma_{ijk}\geq 0, \gamma_{jik}\geq 0, \gamma_{kij}\geq 0,$$
the energy function $\mathcal{E}_{ijk}(u)$ defined on $\mathcal{U}^E_{ijk}$ by (\ref{local Euclid energy function})
could be extended to the following function
\begin{equation}\label{extension function Euclid}
\widetilde{\mathcal{E}}_{ijk}(u)=\int_{u_0}^u\widetilde{\theta}_idu_i+\widetilde{\theta}_jdu_j+\widetilde{\theta}_kdu_k,
\end{equation}
which is a $C^1$-smooth concave function defined on $\mathbb{R}^3$ with
$$\nabla_u \widetilde{\mathcal{E}}_{ijk}=(\widetilde{\theta}_i, \widetilde{\theta}_j, \widetilde{\theta}_k)^T.$$
\end{lemma}

Using Lemma \ref{Euclidean extension lemma}, we can prove the following global rigidity of Euclidean inversive distance circle packings,
which is the Euclidean part of Theorem \ref{main theorem global rigidity for Curvature}.

\begin{theorem}\label{Euclidean rigidity}
Given a triangulated surface $(M, \mathcal{T})$ with inversive distance $I>-1$ and
$\gamma_{ijk}\geq 0, \gamma_{jik}\geq 0, \gamma_{kij}\geq 0$
for any topological triangle $\triangle ijk\in F$. Then for any $\overline{K}\in C(V)$ with $\sum_{i=1}^N\overline{K}_i=2\pi\chi(M)$,
there exists at most one Euclidean inversive distance circle packing metric $r$ up to scaling with $K(r)=\overline{K}$.
\end{theorem}

\proof
By Lemma \ref{Euclidean extension lemma},
the Ricci potential function $\mathcal{E}(u)$ in (\ref{local Euclid energy function}) could be extended from $\mathcal{U}^E$ to the whole space $\mathbb{R}^{N}$ as follows
$$\widetilde{\mathcal{E}}(u)=-\sum_{\triangle ijk\in F}\widetilde{\mathcal{E}}_{ijk}(u)+\int_{u_0}^u\sum_{i=1}^N(2\pi-\overline{K}_i)du_i.$$
As $\widetilde{\mathcal{E}}_{ijk}(u)$ is $C^1$-smooth concave by Lemma \ref{Euclidean extension lemma} and $\int_{u_0}^u\sum_{i=1}^N(2\pi-\overline{K}_i)du_i$
is a well-defined convex function on $\mathbb{R}^N$,
we have $\widetilde{\mathcal{E}}(u)$ is a $C^1$-smooth convex function on $\mathbb{R}^N$.
By Corollary \ref{positivity of Lambda E}, we have $\widetilde{\mathcal{E}}(u)$ is locally strictly convex on $\mathcal{U}^E\cap \{\sum_{i=1}^Nu_i=0\}$.
Furthermore,
\begin{equation*}
\nabla_{u_i}\widetilde{\mathcal{E}}=-\sum_{\triangle ijk\in F}\widetilde{\theta}_i+2\pi-\overline{K}_i=\widetilde{K}_i-\overline{K}_i,
\end{equation*}
where $\widetilde{K}_i=2\pi-\sum_{\triangle ijk\in F}\widetilde{\theta}_i$,
which implies that $r\in \Omega^E$ is a metric with curvature $\overline{K}$
if and only if the corresponding $u\in \mathcal{U}^E$ is a critical point of $\widetilde{\mathcal{E}}$.

If there are two different inversive distance circle packing metrics $\overline{r}_{A}, \overline{r}_{B}\in \Omega^E$ with the same combinatorial
Curvature $\overline{K}$, then
$\overline{u}_A=\ln \overline{r}_{A}\in \mathcal{U}^E$, $\overline{u}_B=\ln \overline{r}_{B}\in \mathcal{U}^E$
are both critical points of the extended Ricci potential $\widetilde{\mathcal{E}}(u)$.
It follows that
$$\nabla \widetilde{\mathcal{E}}(\overline{u}_A)=\nabla \widetilde{\mathcal{E}}(\overline{u}_B)=0.$$
Set
\begin{equation*}
\begin{aligned}
f(t)=&\widetilde{\mathcal{E}}((1-t)\overline{u}_A+t\overline{u}_B)\\
=&\sum_{\triangle ijk\in F}f_{ijk}(t)+\int_{u_0}^{(1-t)\overline{u}_A+t\overline{u}_B}\sum_{i=1}^N(2\pi-\overline{K}_i)du_i,
\end{aligned}
\end{equation*}
where
$$f_{ijk}(t)=-\widetilde{\mathcal{E}}_{ijk}((1-t)\overline{u}_A+t\overline{u}_B).$$
Then $f(t)$ is a $C^1$ convex function on $[0, 1]$ and $f'(0)=f'(1)=0$, which implies that $f'(t)\equiv 0$ on $[0, 1]$.
Note that $\overline{u}_A$ belongs to the open set $\mathcal{U}^E$,
so there exists $\epsilon>0$ such that $(1-t)\overline{u}_A+t\overline{u}_B\in \mathcal{U}^E$ for $t\in [0, \epsilon]$
and $f(t)$ is smooth on $[0, \epsilon]$.

Note that $f(t)$ is $C^1$ convex on $[0, 1]$ and smooth on $[0, \epsilon]$.
$f'(t)\equiv 0$ on $[0, 1]$ implies that $f''(t)\equiv 0$ on $[0, \epsilon]$.
Note that, for $t\in [0, \epsilon]$,
\begin{equation*}
\begin{aligned}
f''(t)=(\overline{u}_A-\overline{u}_B) \Lambda^E  (\overline{u}_A-\overline{u}_B)^T,
\end{aligned}
\end{equation*}
where $\Lambda^E=-\sum_{\triangle ijk\in F}\Lambda^E_{ijk}$.
By Corollary \ref{positivity of Lambda E}, we have $\overline{u}_A-\overline{u}_B=c(1,\cdots,1)$ for some constant $c\in \mathbb{R}$, which
implies that $\overline{r}_A=e^{c/2}\overline{r}_B$. So there exists at most one Euclidean inversive distance circle packing metric
with combinatorial curvature $\overline{K}$ up to scaling.
\qed

\begin{remark}
The proof of Theorem \ref{Euclidean rigidity} is based on a variational principle, which was introduce by Colin de Verdiere \cite{DV}.
Guo \cite{Guo} used the variational principle to study the infinitesimal rigidity of inversive distance
circle packing metrics for nonnegative inversive distances. Bobenko, Pinkall and Springborn \cite{BPS} introduced a
method to extend a local convex function on a nonconvex domain to a convex function and solved
affirmably a conjecture of Luo \cite{L1} on the global rigidity of piecewise linear metrics. Based on the extension method,
Luo \cite{L3} proved the global rigidity of inversive distance circle packing metrics for nonnegative inversive distance
using the variational principle.
\end{remark}

\subsection{Rigidity of combinatorial $\alpha$-curvature in Euclidean background geometry}

As noted in \cite{GX4}, the classical definition of combinatorial curvature $K_i$ with Euclidean background geometry
in (\ref{classical Gauss curv}) has two disadvantages. The first is that the classical combinatorial curvature is scaling
invariant, i.e. $K_i(\lambda r)=K_i(r)$ for any $\lambda>0$; The second is that, as the triangulated surfaces approximate a smooth surface,
the classical combinatorial curvature $K_i$ could not approximate the smooth Gauss curvature, as we obviously have $K_i$ tends zero.
Motivated by the observations, Ge and the author introduced a new combinatorial curvature for triangulated surfaces with Thurston's circle
packing metrics in \cite{GX4,GX3,GX5}. Ge and Jiang \cite{GJ4} and Ge and the author \cite{GX6} further generalized the curvature to inversive distance circle packing metrics.
Set
\begin{equation}\label{s_i}
\begin{aligned}
s_i(r)=\left\{
         \begin{array}{ll}
           r_i, & \hbox{Euclidean background geometry} \\
           \tanh\frac{r_i}{2}, & \hbox{hyperbolic background geometry}
         \end{array}
       \right.
\end{aligned}.
\end{equation}
We have the following definition of combinatorial $\alpha$-curvature on triangulated surfaces
with inversive distance circle packing metrics.

\begin{definition}\label{definition of comb curv with inversive dist}
Given a triangulated surface $(M, \mathcal{T})$ with inversive distance $I>-1$
and an inversive distance circle packing metric $r\in \Omega$, the combinatorial $\alpha$-curvature at the vertex $i$ is defined
to be
\begin{equation}\label{definition of Gauss curv}
R_{\alpha,i}=\frac{K_i}{s_i^\alpha},
\end{equation}
where $\alpha\in \mathbb{R}$ is a constant,
$K_i$ is the classical combinatorial curvature at $i$ given by (\ref{classical Gauss curv}) and $s_i$ is given by (\ref{s_i}).
\end{definition}

Specially, if $\alpha=0$, then $R_{\alpha,i}=K_i$.
As the inversive distance generalizes Thurston's intersection angle, the Definition \ref{definition of comb curv with inversive dist}
of combinatorial $\alpha$-curvature naturally generalizes the definition of combinatorial curvature in \cite{GX4,GX3,GX5}.

For the $\alpha$-curvature $R_{\alpha,i}$, we have the following global rigidity of Euclidean inversive distance circle packing metrics
for inversive distance in $(-1, +\infty)$, which is the Euclidean part of Theorem \ref{main theorem global rigidity for alpha curvature}.

\begin{theorem}\label{Euclidean alpha global rigidity}
Given a closed triangulated surface $(M, \mathcal{T})$ with inversive distance $I>-1$ and
$\gamma_{ijk}\geq 0, \gamma_{jik}\geq 0, \gamma_{kij}\geq 0$
for any topological triangle $\triangle ijk\in F$. $\overline{R}$
is a given function defined on the vertices  of  $(M, \mathcal{T})$.
If $\alpha\overline{R}\equiv0$, there exists at most one Euclidean inversive
distance circle packing metric $\overline{r}\in \Omega^E$
with $\alpha$-curvature $\overline{R}$ up to scaling. If $\alpha\overline{R}\leq0$ and $\alpha\overline{R}\not\equiv0$,
there exists at most one Euclidean inversive distance circle packing metric $\overline{r}\in \Omega^E$ with $\alpha$-curvature $\overline{R}$.
\end{theorem}

As the proof of Theorem \ref{Euclidean alpha global rigidity} is almost parallel to that of Theorem \ref{Euclidean rigidity} using the energy function
$$\widetilde{\mathcal{E}}_\alpha(u)=-\sum_{ \triangle ijk\in F}\widetilde{\mathcal{E}}_{ijk}(u)+\int_{u_0}^u\sum_{i=1}^N(2\pi-\overline{R}_ir_i^\alpha)du_i,$$
we omit the details of the proof.

\section{Hyperbolic inversive distance circle packing metrics}\label{section 4}
\subsection{Admissible space of hyperbolic inversive distance circle packing metrics for a single triangle}
In this subsection, we investigate the admissible space of hyperbolic inversive distance circle packings
for a single topological triangle $\triangle ijk\in F$ with inversive distance $I_i, I_j, I_k\in (-1, +\infty)$
and
\begin{equation}\label{inversive condition}
\gamma_{ijk}\geq 0, \gamma_{jik}\geq 0, \gamma_{kij}\geq 0.
\end{equation}
Suppose $\triangle ijk$ is a topological triangle in $F$.
In the hyperbolic background geometry, the length $l_{i}$ of the edge $\{jk\}$ is defined by
\begin{equation}\label{hyperbolic length}
l_{i}=\cosh^{-1}(\cosh r_j\cosh r_k+I_i\sinh r_j\sinh r_k),
\end{equation}
where $I_i$ is the hyperbolic inversive distance between the two circles attached to the vertices $j$ and $k$.
In order that the edge lengths $l_i, l_j, l_k$ satisfy the triangle inequalities, there are some restrictions on the
radius vectors. So we first study the triangle inequalities for the hyperbolic background geometry.
To simplify the notations, we use the following simplification
$$C_i=\cosh r_i, S_i=\sinh r_i,$$
when there is no confusion. We have the following lemma on the hyperbolic triangle inequalities.

\begin{lemma}\label{hyperbolic triangle inequality}
Suppose $(M, \mathcal{T}, I)$ is a weighted triangulated surface with hyperbolic inversive distance $I>-1$ and
$\triangle ijk$ is a topological triangle in $F$.
Suppose $l_i, l_j, l_k$ are the edge lengths defined by the hyperbolic inversive distance $I_i, I_j, I_k$ using the
radius $r_i, r_j, r_k$ by (\ref{hyperbolic length}), then the triangle inequalities are satisfied if and only if
\begin{equation}
\begin{aligned}
&2S_i^2S_j^2S_k^2(1+I_iI_jI_k)+S_i^2S_j^2(1-I_k^2)+S_i^2S_k^2(1-I_j^2)+S_j^2S_k^2(1-I_i^2)\\
&+2C_jC_kS_i^2S_jS_k\gamma_{ijk}+2C_iC_kS_iS_j^2S_k\gamma_{jik}+2C_iC_jS_iS_jS_k^2\gamma_{kij}> 0.
\end{aligned}
\end{equation}
\end{lemma}
\proof
In order that
$l_i+l_j>l_k, l_i+l_k>l_j, l_j+l_k>l_i, $
we just need
$$\sinh\frac{l_i+l_j-l_k}{2}>0, \sinh\frac{l_i+l_k-l_j}{2}>0, \sinh\frac{l_j+l_k-l_i}{2}>0.$$
Note that $l_i>0, l_j>0, l_k>0$, this is equivalent to
$$\sinh\frac{l_i+l_j+l_k}{2}\sinh\frac{l_i+l_j-l_k}{2}\sinh\frac{l_i+l_k-l_j}{2}\sinh\frac{l_j+l_k-l_i}{2}>0.$$
By direct calculations, we have
\begin{equation*}
\begin{aligned}
&4\sinh\frac{l_i+l_j+l_k}{2}\sinh\frac{l_i+l_j-l_k}{2}\sinh\frac{l_i+l_k-l_j}{2}\sinh\frac{l_j+l_k-l_i}{2}\\
=&(\cosh (l_i+l_j)-\cosh l_k)(\cosh l_k-\cosh (l_i-l_j))\\
=&(\cosh^2 l_i-1)(\cosh l_j^2-1)-(\cosh l_i\cosh l_j-\cosh l_k)^2\\
=&(2C_i^2C_j^2C_k^2-C_i^2C_j^2-C_i^2C_k^2-C_j^2C_k^2+1)-(S_i^2S_j^2I_k^2+S_i^2S_k^2I_j^2+S_j^2S_k^2I_i^2)\\
 &+2C_jC_kS_i^2S_jS_kI_i+2C_iC_kS_iS_j^2S_kI_j+2C_iC_jS_iS_jS_k^2I_k\\
 &+2C_iC_jS_iS_jS_k^2I_iI_j+2C_iC_kS_iS_j^2S_kI_iI_k+2C_jC_kS_i^2S_jS_kI_jI_k+2S_i^2S_j^2S_k^2I_iI_jI_k,
\end{aligned}
\end{equation*}
where the definition of edge length (\ref{hyperbolic length}) is
used in the last line.
Note that
$$C_i^2=\cosh^2r_i=\sinh^2r_i+1=S_i^2+1,$$
we have
\begin{equation*}
\begin{aligned}
&4\sinh\frac{l_i+l_j+l_k}{2}\sinh\frac{l_i+l_j-l_k}{2}\sinh\frac{l_i+l_k-l_j}{2}\sinh\frac{l_j+l_k-l_i}{2}\\
=&2S_i^2S_j^2S_k^2(1+I_iI_jI_k)+S_i^2S_j^2(1-I_k^2)+S_i^2S_k^2(1-I_j^2)+S_j^2S_k^2(1-I_i^2)\\
&+2C_jC_kS_i^2S_jS_k(I_i+I_jI_k)+2C_iC_kS_iS_j^2S_k(I_j+I_iI_k)+2C_iC_jS_iS_jS_k^2(I_k+I_iI_j).
\end{aligned}
\end{equation*}
This completes the proof of the lemma.
\qed

Denote the admissible space of hyperbolic inversive distance circle packing metrics for a triangle $\triangle ijk\in F$
as $\Omega^{H}_{ijk}$, i.e.
$$\Omega^{H}_{ijk}:=\{(r_i,r_j,r_k)\in \mathbb{R}^3_{>0}|l_i+l_j>l_k, l_i+l_k>l_j, l_j+l_k>l_i \}.$$
By Lemma \ref{hyperbolic triangle inequality}, we have the following direct corollary, which was obtained by Zhou \cite{Z}.

\begin{corollary}
Suppose $\triangle ijk$ is a topological triangle in $F$ with hyperbolic inversive distance $I_i, I_j, I_k\in (-1, 1]$ and
$\gamma_{ijk}\geq 0, \gamma_{jik}\geq 0, \gamma_{kij}\geq 0,$
then
$\Omega^H_{ijk}=\mathbb{R}^{3}_{>0}$, i.e.
the triangle inequalities are satisfied for all radius vectors in $\mathbb{R}^3_{>0}$.
\end{corollary}

Specially, if $I_i=\cos \Phi_{i}, I_j=\cos \Phi_{j}, I_k=\cos \Phi_{k}$
with $\Phi_{i}, \Phi_{j}, \Phi_{k}\in [0, \frac{\pi}{2}]$, the triangle inequalities are satisfied for all radius vectors,
which was obtained by Thurston in \cite{T1}.

By Lemma \ref{hyperbolic triangle inequality}, we can also get the following useful result.

\begin{corollary}\label{hyperbolic triangle inequality is valid for large s}
Suppose $\triangle ijk$ is a topological triangle in $F$ with hyperbolic inversive distance $I>-1$ and
$\gamma_{ijk}\geq 0, \gamma_{jik}\geq 0, \gamma_{kij}\geq 0.$
Suppose the edge lengths $l_i, l_j, l_k$ are generated by the radius vector $(s,s,s)$ with $s\in \mathbb{R}_{>0}$.
If $s\in \mathbb{R}_{>0}$ satisfies
\begin{equation}\label{hyperbolic triangle inequ for s large}
\begin{aligned}
\sinh^2 s>\frac{I_i^2+I_j^2+I_k^2-3}{2(1+I_i)(1+I_j)(1+I_j)}.
\end{aligned}
\end{equation}
we have $(s,s,s)\in \Omega^H_{ijk}$.
\end{corollary}

\proof
By Lemma \ref{hyperbolic triangle inequality}, for $s>0$, $(s,s,s)\in \Omega^H_{ijk}$ if and only if
\begin{equation*}
\begin{aligned}
2\cosh^2s(\gamma_{ijk}+\gamma_{jik}+\gamma_{kij})+2\sinh^2 s(1+I_iI_jI_k)+3-I_i^2-I_j^2-I_k^2>0.
\end{aligned}
\end{equation*}
By $\gamma_{ijk}\geq 0$, $\gamma_{jik}\geq 0$, $\gamma_{kij}\geq 0$,
we have $\gamma_{ijk}+\gamma_{jik}+\gamma_{kij}\geq 0$. Then
\begin{equation*}
\begin{aligned}
&2\cosh^2s(\gamma_{ijk}+\gamma_{jik}+\gamma_{kij})+2\sinh^2 s(1+I_iI_jI_k)+3-I_i^2-I_j^2-I_k^2\\
\geq&2\sinh^2 s(1+I_iI_jI_k+\gamma_{ijk}+\gamma_{jik}+\gamma_{kij})+3-I_i^2-I_j^2-I_k^2\\
=&2\sinh^2 s(1+I_i)(1+I_j)(1+I_j)+3-I_i^2-I_j^2-I_k^2.
\end{aligned}
\end{equation*}
Note that $I_{i}, I_j, I_k\in (-1, +\infty)$, to ensure the triangle inequalities, we just need
\begin{equation*}
\begin{aligned}
\sinh^2 s>\frac{I_i^2+I_j^2+I_k^2-3}{2(1+I_i)(1+I_j)(1+I_j)}.
\end{aligned}
\end{equation*}
 \qed

Guo \cite{Guo} obtained a result similar to Corollary \ref{hyperbolic triangle inequality is valid for large s} for $I\geq 0$.

By Lemma \ref{hyperbolic triangle inequality}, $\Omega^H_{ijk}\neq \mathbb{R}^3_{>0}$ for general $I_i, I_j, I_k\in (-1, +\infty)$.
Furthermore, $\Omega^H_{ijk}$ is not convex. Similar to the case of Euclidean background geometry,
we have the following lemma on the structure of $\Omega^{H}_{ijk}$.
\begin{lemma}\label{hyperbolic simply connect}
Suppose $\triangle ijk$ is a topological triangle in $F$ with hyperbolic inversive distance $I>-1$ and
$\gamma_{ijk}\geq 0, \gamma_{jik}\geq 0, \gamma_{kij}\geq 0,$
then the admissible space $\Omega^{H}_{ijk}$ is simply connected.
Furthermore, for each connected component $V$ of $\mathbb{R}^3_{>0}\setminus \Omega^H_{ijk}$,
the intersection $V\cap \overline{\Omega}^H_{ijk}$ is a connected component of $\overline{\Omega}^H_{ijk}\setminus \Omega^H_{ijk}$,
on which $\theta_i$ is a constant function.
\end{lemma}
\proof
Define the map
\begin{equation*}
\begin{aligned}
F: \mathbb{R}^3_{>0}&\rightarrow \mathbb{R}^3_{>0}\\
(r_i,r_j,r_k)&\mapsto (F_i, F_j, F_k)
\end{aligned}
\end{equation*}
where
\begin{equation*}
\begin{aligned}
F_i&=\cosh r_j\cosh r_k+I_i\sinh r_j\sinh r_k,\\
F_j&=\cosh r_i\cosh r_k+I_j\sinh r_i\sinh r_k,\\
F_k&=\cosh r_i\cosh r_j+I_k\sinh r_i\sinh r_j.
\end{aligned}
\end{equation*}
By direct calculations, we have
\begin{equation*}
\begin{aligned}
\frac{\partial (F_i, F_j, F_k)}{\partial (r_i, r_j, r_k)}
=
\left(
  \begin{array}{ccc}
    0 & S_jC_k+I_iC_jS_k & C_jS_k+I_iS_jC_k \\
    S_iC_k+I_jC_iS_k & 0 & C_iS_k+I_jS_iC_k \\
    S_iC_j+I_kC_iS_j & C_iS_j+I_kS_iC_j & 0 \\
  \end{array}
\right)
\end{aligned}
\end{equation*}
and
\begin{equation*}
\begin{aligned}
\left|\frac{\partial (F_i, F_j, F_k)}{\partial (r_i, r_j, r_k)}\right|
=&2C_iC_jC_kS_iS_jS_k(1+I_iI_jI_k)+\gamma_{kij}C_kS_k(C_i^2S_j^2+C_j^2S_i^2)\\
 &+\gamma_{jik}C_jS_j(C_k^2S_i^2+C_i^2S_k^2)+\gamma_{ijk}C_iS_i(C_k^2S_j^2+C_j^2S_k^2).
\end{aligned}
\end{equation*}
By $I>-1$ and
$\gamma_{ijk}\geq 0, \gamma_{jik}\geq 0, \gamma_{kij}\geq 0,$
we have
\begin{equation*}
\begin{aligned}
\left|\frac{\partial (F_i, F_j, F_k)}{\partial (r_i, r_j, r_k)}\right|
\geq &2C_iC_jC_kS_iS_jS_k(1+I_iI_jI_k+\gamma_{ijk}+\gamma_{jik}+\gamma_{kij})\\
=&2C_iC_jC_kS_iS_jS_k(1+I_i)(1+I_j)(1+I_k)>0,
\end{aligned}
\end{equation*}
which implies that $F$ is globally injective. In fact, if there are two different $r=(r_i, r_j, r_k)$ and $r'=(r'_i, r'_j, r'_k)$ satisfying $F(r)=F(r')$, then we have
\begin{equation*}
\begin{aligned}
0=F(r)-F(r')=\frac{\partial (F_i, F_j, F_k)}{\partial (r_i, r_j, r_k)}\Huge{|}_{r+\theta(r-r')}\cdot(r-r')^T, 0<\theta<1,
\end{aligned}
\end{equation*}
which implies $r=r'$ by the nondegeneracy of $\frac{\partial (F_i, F_j, F_k)}{\partial (r_i, r_j, r_k)}$ on $\mathbb{R}^3_{>0}$.
So the map $F$ is injective on $\mathbb{R}^3_{>0}$.

Note that $F$ has the following property
\begin{equation*}
\begin{aligned}
0<(1+I_i)\sinh r_j\sinh r_k\leq F_i\leq (1+|I_i|)\cosh (r_i+r_j),
\end{aligned}
\end{equation*}
which implies that $F$ is a proper map.
By the invariance of domain theorem, we have
$F: \mathbb{R}^3_{>0}\rightarrow F(\mathbb{R}^3_{>0})$ is a diffeomorphism.

Define
\begin{equation*}
\begin{aligned}
G: \mathbb{R}^3_{>0}&\rightarrow \mathbb{R}^3_{>0}\\
(l_i,l_j,l_k)&\mapsto (\cosh l_i, \cosh l_j, \cosh l_k),
\end{aligned}
\end{equation*}
then $G: \mathbb{R}^3_{>0}\rightarrow G(\mathbb{R}^3_{>0})$ is a diffeomorphis and
$H=G^{-1}\circ F$ is the map defining the edge length by the inversive distance
which maps $(r_i,r_j, r_k)$ to $(l_i, l_j, l_k)$.

Set
$$\mathcal{L}=\{(l_i,l_j,l_k)|l_i+l_j>l_k, l_i+l_k>l_j, l_j+l_k>l_i\},$$
then $\Omega^H_{ijk}=H^{-1}(H(\mathbb{R}^3_{>0})\cap \mathcal{L})$. To prove that $\Omega^H_{ijk}$
is simply connected, we just need to prove that $H(\mathbb{R}^3_{>0})\cap \mathcal{L}$ is simply connected.

Note that $\mathcal{L}$ is a cone in $\mathbb{R}^3_{>0}$ bounded by three planes
\begin{equation*}
\begin{aligned}
L_i=&\{(l_{i}, l_j, l_k)\in \mathbb{R}^3_{>0}|l_i=l_j+l_k\},\\
L_j=&\{(l_{i}, l_j, l_k)\in \mathbb{R}^3_{>0}|l_j=l_i+l_k\},\\
L_k=&\{(l_{i}, l_j, l_k)\in \mathbb{R}^3_{>0}|l_k=l_i+l_j\}.\\
\end{aligned}
\end{equation*}
By the fact that $H$ is a diffeomorphism between $\mathbb{R}^3_{>0}$ and $H(\mathbb{R}^3_{>0})$,
$H(\mathbb{R}^3_{>0})$ is the set bounded by three surfaces
\begin{equation*}
\begin{aligned}
\Sigma_i=&\{(l_{i}, l_j, l_k)\in \mathbb{R}^3_{>0}|\cosh l_i=\cosh l_j\cosh l_k+I_i\sinh l_j\sinh l_k\},\\
\Sigma_j=&\{(l_{i}, l_j, l_k)\in \mathbb{R}^3_{>0}|\cosh l_j=\cosh l_i\cosh l_k+I_j\sinh l_i\sinh l_k\},\\
\Sigma_k=&\{(l_{i}, l_j, l_k)\in \mathbb{R}^3_{>0}|\cosh l_k=\cosh l_i\cosh l_j+I_k\sinh l_i\sinh l_j\}.
\end{aligned}
\end{equation*}
In fact, if $r_i=0$, then $l_j=r_k, l_k=r_j$ and
$\cosh l_i=\cosh r_j\cosh r_k+I_i\sinh r_j\sinh r_k=\cosh l_j\cosh l_k+I_i\sinh l_j\sinh l_k.$
$\Sigma_i$ is in fact the image of $r_i=0$ under $H$.
By the diffeomorphism of $H$, $\Sigma_i$, $\Sigma_j$, $\Sigma_k$ are mutually disjoint.
Furthermore, if $I_i\in (-1, 1]$, we have
$\cosh(l_j-l_k)<\cosh l_i\leq\cosh (l_j+l_k)$ on $\Sigma_i$. And if $I_i\in (1, +\infty)$, we have $\cosh l_i>\cosh(l_j+l_k)$ on $\Sigma_i$.
This implies that $\Sigma_i\subset \overline{\mathcal{L}}$ if $I_i\in (-1, 1]$ and $\Sigma_i\cap \mathcal{L}=\emptyset$ if $I_i\in (1, +\infty)$.
Similar results hold for $\Sigma_j$ and $\Sigma_k$.
To prove that $H(\mathbb{R}^3_{>0})\cap \mathcal{L}$ is simply connected, we just need to consider the following cases by the symmetry between $i,j,k$.

If $I_i, I_j, I_k\in (-1, 1]$, $H(\mathbb{R}^3_{>0})\cap \mathcal{L}$ is bounded by $\Sigma_i, \Sigma_j, \Sigma_k$ and
$H(\mathbb{R}^3_{>0})\cap \mathcal{L}=H(\mathbb{R}^3_{>0})$.

If $I_i, I_j\in (-1, 1]$ and $I_k\in(1,+\infty)$, $H(\mathbb{R}^3_{>0})\cap \mathcal{L}$ is bounded by $\Sigma_i, \Sigma_j$ and $L_k$.

If $I_i\in (-1, 1]$ and $I_j, I_k\in(1,+\infty)$, $H(\mathbb{R}^3_{>0})\cap \mathcal{L}$ is bounded by $\Sigma_i$, $L_j$ and $L_k$.

If $I_i, I_j, I_k\in(1,+\infty)$, $H(\mathbb{R}^3_{>0})\cap \mathcal{L}$ is bounded by $L_i$, $L_j$ and $L_k$.
In this case, $H(\mathbb{R}^3_{>0})\cap \mathcal{L}=\mathcal{L}$.

For any case, $H(\mathbb{R}^3_{>0})\cap \mathcal{L}$ is a simply connected subset of $\mathbb{R}^3_{>0}$.
By the fact that $H$ is a diffeomorphism between $\mathbb{R}^3_{>0}$ and $H(\mathbb{R}^3_{>0})$,
we have the admissible space $\Omega^H_{ijk}=H^{-1}(H(\mathbb{R}^3_{>0})\cap \mathcal{L})$ is simply connected.

By the analysis above,
if $H(\mathbb{R}^3_{>0})\subset \mathcal{L}$, then $\Omega^H_{ijk}=H^{-1}(H(\mathbb{R}^3_{>0})\cap\mathcal{L})=\mathbb{R}^3_{>0}$.
If $H(\mathbb{R}^3_{>0})\setminus\mathcal{L}\neq \emptyset$, then $\Omega^H_{ijk}$ is a proper subset of $\mathbb{R}^3_{>0}$.
If $I_i>1$, the boundary component
$\Sigma_i=\{(l_{i}, l_j, l_k)\in \mathbb{R}^3_{>0}|\cosh l_i=\cosh l_j\cosh l_k+I_i\sinh l_j\sinh l_k\}$
is out of the set $\mathcal{L}$.
By the fact that $\Omega^H_{ijk}=H^{-1}(H(\mathbb{R}^3_{>0})\cap\mathcal{L})$ and
$H:\mathbb{R}^3_{>0}\rightarrow H(\mathbb{R}^3_{>0})$ is a diffeomorphism,
we have $H^{-1}(L_i)$ is a connected boundary component of $\Omega^H_{ijk}$, on which $\theta_i=\pi, \theta_j=\theta_k=0$.
This completes the proof of the lemma.
\qed

\begin{corollary}\label{hyperbolic triangle extension lemma}
For a topological triangle $\triangle ijk\in F$ with inversive distance $I>-1$ and
$\gamma_{ijk}\geq 0, \gamma_{jik}\geq 0, \gamma_{kij}\geq 0,$
the functions $\theta_i, \theta_j, \theta_k$ defined on $\Omega^H_{ijk}$
could be continuously extended by constant to
$\widetilde{\theta}_i, \widetilde{\theta}_j, \widetilde{\theta}_k$  defined on $\mathbb{R}^3_{>0}$.
\end{corollary}

\subsection{Infinitesimal rigidity of hyperbolic inversive distance circle packings}

Set $u_i=\ln \tanh \frac{r_i}{2}$, then we have $\mathcal{U}^H_{ijk}:=u (\Omega^H_{ijk})$ is a simply connected subset of $\mathbb{R}^3_{>0}$.
If $(r_i, r_j, r_k)\in \Omega^H_{ijk}$, $l_i, l_j,l_k$ form a hyperbolic triangle.
Denote the inner angle at the vertex $i$ as $\theta_i$. We have the following lemma.

\begin{lemma}\label{hyperbolic symmetry lemma}
For any triangle $\triangle ijk\in F$, we have
\begin{equation}
\frac{\partial \theta_i}{\partial u_j}=\frac{\partial \theta_j}{\partial u_i}=\frac{1}{A\sinh^2 l_{k}}[C_kS_i^2S_j^2(1-I_k^2)+C_iS_iS_j^2S_k\gamma_{jik}+C_jS^2_iS_jS_k\gamma_{ijk}]
\end{equation}
on $\mathcal{U}^H_{ijk}$, where $A=\sinh l_j\sinh l_k\sin \theta_i$.
\end{lemma}
\proof
By cosine law, we have
$\cosh l_i=\cosh l_j\cosh l_k-\sinh l_j\sinh l_k\cos\theta_i.$
Taking the derivative with respect to $l_i$ gives
$$\frac{\partial \theta_i}{\partial l_i}=\frac{\sinh l_i}{A},$$
where $A=\sinh l_j\sinh l_k\sin \theta_i$.
Similarly, taking the derivative with respect to $l_j$ and $l_k$ and using the cosine law again, we have
\begin{equation*}
\begin{aligned}
\frac{\partial \theta_i}{\partial l_j}=\frac{-\sinh l_i\cos\theta_k}{A},  \frac{\partial \theta_i}{\partial l_k}=\frac{-\sinh l_i\cos\theta_j}{A}.
\end{aligned}
\end{equation*}
By the definition of edge length $l_i, l_j$ and $l_k$, we have
$$\frac{\partial l_i}{\partial r_j}=\frac{\sinh r_j\cosh r_k+I_i\cosh r_j\sinh r_k}{\sinh l_i},
\frac{\partial l_j}{\partial r_j}=0,
\frac{\partial l_k}{\partial r_j}=\frac{\sinh r_j\cosh r_i+I_k\cosh r_j\sinh r_i}{\sinh l_k}.$$
Then
\begin{equation*}
\begin{aligned}
A\frac{\partial \theta_i}{\partial u_j}
=&A\sinh r_j\frac{\partial \theta_i}{\partial r_j}\\
=&A\sinh r_j(\frac{\partial \theta_i}{\partial l_i}\frac{\partial l_i}{\partial r_j}+\frac{\partial \theta_i}{\partial l_k}\frac{\partial l_k}{\partial r_j})\\
=&\sinh r_j(\sinh r_j\cosh r_k+I_i\cosh r_j\sinh r_k)\\
 &-\frac{1}{\sinh l_k}\sinh r_j\sinh l_i\cos\theta_j(\sinh r_j\cosh r_i+I_k\cosh r_j\sinh r_i),
\end{aligned}
\end{equation*}
which implies that
\begin{equation*}
\begin{aligned}
\sinh^2l_kA\frac{\partial \theta_i}{\partial u_j}
=&(\cosh^2l_k-1)\sinh r_j(\sinh r_j\cosh r_k+I_i\cosh r_j\sinh r_k)\\
 &+(\cosh l_j-\cosh l_i\cosh l_k)\sinh r_j(\sinh r_j\cosh r_i+I_k\cosh r_j\sinh r_i)
\end{aligned}
\end{equation*}
Note that
\begin{equation*}
\begin{aligned}
\sinh r_j(\sinh r_j\cosh r_k+I_i\cosh r_j\sinh r_k)=&\cosh r_j\cosh l_i-\cosh r_k,\\
\sinh r_j(\sinh r_j\cosh r_i+I_k\cosh r_j\sinh r_i)=&\cosh r_j\cosh l_k-\cosh r_i.
\end{aligned}
\end{equation*}
Using the definition of of edge lengths $l_i, l_j$ and $l_k$, by direct calculations, we have
\begin{equation*}
\begin{aligned}
\frac{\partial \theta_i}{\partial u_j}
=\frac{1}{A\sinh^2 l_{k}}[C_kS_i^2S_j^2(1-I_k^2)+C_iS_iS_j^2S_k\gamma_{jik}+C_jS^2_iS_jS_k\gamma_{ijk}],
\end{aligned}
\end{equation*}
which implies also $\frac{\partial \theta_i}{\partial u_j}=\frac{\partial \theta_j}{\partial u_i}$.
\qed

\begin{remark}\label{positivity of hyperbolic derivative}
For $I_{i}, I_j, I_k\in (-1, 1]$ and $\gamma_{ijk}\geq 0$, $\gamma_{jik}\geq 0$, $\gamma_{kij}\geq 0$,
by Lemma \ref{hyperbolic symmetry lemma},
we have $\frac{\partial \theta_i}{\partial u_j}\geq 0$, and $\frac{\partial \theta_i}{\partial u_j}=0$ if and only if
$I_k=1$ and $I_i+I_j=0$. Especially, if $I_i=\cos\Phi_{i}, I_j=\cos\Phi_{j}, I_{k}=\cos\Phi_{k}$
with $\Phi_{i},\Phi_{j},\Phi_{k}\in [0,\frac{\pi}{2}]$, we have
$\frac{\partial \theta_i}{\partial u_j}\geq 0$, and $\frac{\partial \theta_i}{\partial u_j}=0$ if and only if
$\Phi_{k}=0$ and $\Phi_{i}=\Phi_{j}=\frac{\pi}{2}$.
\end{remark}

Lemma \ref{hyperbolic symmetry lemma} shows that the matrix
\begin{equation*}
\begin{aligned}
\Lambda^H_{ijk}=\frac{\partial (\theta_i, \theta_j, \theta_k)}{\partial (u_i, u_j, u_k)}
=\left(
   \begin{array}{ccc}
     \frac{\partial \theta_i}{\partial u_i} & \frac{\partial \theta_i}{\partial u_j} & \frac{\partial \theta_i}{\partial u_k} \\
     \frac{\partial \theta_j}{\partial u_i} & \frac{\partial \theta_j}{\partial u_j} & \frac{\partial \theta_j}{\partial u_k} \\
     \frac{\partial \theta_k}{\partial u_i} & \frac{\partial \theta_k}{\partial u_j} & \frac{\partial \theta_k}{\partial u_k} \\
   \end{array}
 \right)
\end{aligned}
\end{equation*}
is symmetric on $\mathcal{U}^H_{ijk}$. Similar to the case of Euclidean background geometry, we have the following lemma
for the matrix $\Lambda^H_{ijk}$.

\begin{lemma}\label{hyperbolic negative definite}
In the hyperbolic background geometry, for any triangle $\triangle ijk\in F$ with $I_i, I_j, I_k>-1$ and
$\gamma_{ijk}\geq 0$, $\gamma_{jik}\geq 0$, $\gamma_{kij}\geq 0$,
the matrix $\Lambda^H_{ijk}$ is negative definite on $\mathcal{U}^H_{ijk}$.
\end{lemma}
\proof
The proof is parallel to that of Lemma 12 in \cite{Guo} with some modifications.
By the proof of Lemma \ref{hyperbolic symmetry lemma}, we have
\begin{equation}\label{hyperbolic diff theta}
\begin{aligned}
\left(
  \begin{array}{c}
    d\theta_i \\
    d\theta_j \\
    d\theta_k \\
  \end{array}
\right)
=&-\frac{1}{A}
 \left(
   \begin{array}{ccc}
     \sinh l_i & 0 & 0 \\
     0 & \sinh l_j & 0 \\
     0 & 0 & \sinh l_k \\
   \end{array}
 \right)
 \left(
   \begin{array}{ccc}
     -1 & \cos \theta_k & \cos\theta_j \\
     \cos \theta_k & -1 & \cos \theta_i \\
     \cos \theta_j & \cos \theta_i & -1 \\
   \end{array}
 \right)\\
 &\times
  \left(
   \begin{array}{ccc}
     \frac{1}{\sinh l_i} & 0 & 0 \\
     0 & \frac{1}{\sinh l_j} & 0 \\
     0 & 0 & \frac{1}{\sinh l_k} \\
   \end{array}
   \right)
   \left(
     \begin{array}{ccc}
       0 & R_{ijk} & R_{ikj} \\
       R_{jik} & 0 & R_{jki} \\
       R_{kij} & R_{kji} & 0 \\
     \end{array}
   \right)\\
   &\times
   \left(
   \begin{array}{ccc}
     \sinh r_i & 0 & 0 \\
     0 & \sinh r_j & 0 \\
     0 & 0 & \sinh r_k \\
   \end{array}
 \right)
 \left(
  \begin{array}{c}
    du_i \\
    du_j \\
    du_k \\
  \end{array}
\right),
\end{aligned}
\end{equation}
where
$$A=\sinh l_i\sinh l_j\sin \theta_k, R_{ijk}=\sinh r_j\cosh r_k+I_i\cosh r_j\sinh r_k.$$
Write the equation (\ref{hyperbolic diff theta}) as
\begin{equation}\label{matrix J}
\begin{aligned}
\left(
  \begin{array}{c}
    d\theta_i \\
    d\theta_j \\
    d\theta_k \\
  \end{array}
\right)
=-\frac{1}{A}\mathcal{J} \left(
  \begin{array}{c}
    du_i \\
    du_j \\
    du_k \\
  \end{array}
\right)
\end{aligned}
\end{equation}
and denote the second and fourth matrix in the product of the right hand side of (\ref{hyperbolic diff theta})
as $\Theta$ and $\mathcal{R}$ respectively.
Then $\Lambda^H_{ijk}$ is negative definite is equivalent to $\mathcal{J}$ is positive definite.

We first prove that $\det\mathcal{J}$ is positive. To prove this, we just need to prove that
$\det(\Theta)$ and $\det\mathcal{R}$ are positive. By direct calculations, we have
\begin{equation*}
\begin{aligned}
\det \Theta
=&-1+\cos\theta_i^2+\cos\theta_j^2+\cos\theta_k^2+2\cos \theta_i\cos \theta_j\cos \theta_k\\
=&4\cos\frac{\theta_i+\theta_j-\theta_k}{2}\cos\frac{\theta_i-\theta_j+\theta_k}{2}
 \cos\frac{\theta_i+\theta_j+\theta_k}{2}\cos\frac{\theta_i-\theta_j-\theta_k}{2}.
\end{aligned}
\end{equation*}
By the Gauss-Bonnet formula for hyperbolic triangles, we have
$$\theta_i+\theta_j+\theta_k=\pi-Area(\triangle ijk),$$
which implies $\frac{\theta_i+\theta_j+\theta_k}{2}, \frac{\theta_i+\theta_j-\theta_k}{2},
\frac{\theta_i-\theta_j+\theta_k}{2}, \frac{\theta_i-\theta_j-\theta_k}{2}\in (-\frac{\pi}{2}, \frac{\pi}{2})$.
Then we have $\det \Theta>0$.

By direct calculations, we have
\begin{equation*}
\begin{aligned}
\det \mathcal{R}
=&R_{ijk}R_{jki}R_{kij}+R_{ikj}R_{jik}R_{kji}\\
=&2C_iC_jC_kS_iS_jS_k(1+I_iI_jI_k)+C_kS_k(I_k+I_iI_j)(C_i^2S_j^2+C_j^2S_i^2)\\
 &+C_jS_j(I_j+I_iI_k)(C_k^2S_i^2+C_i^2S_k^2)+C_iS_i(I_i+I_jI_k)(C_k^2S_j^2+C_j^2S_k^2)\\
\geq &2C_iC_jC_kS_iS_jS_k(1+I_iI_jI_k+I_k+I_iI_j+I_j+I_iI_k+I_i+I_jI_k)\\
=&2C_iC_jC_kS_iS_jS_k(1+I_i)(1+I_j)(1+I_k)>0,
\end{aligned}
\end{equation*}
where the conditions $I_{i}, I_j, I_k\in (-1, +\infty)$ and $\gamma_{ijk}\geq 0$, $\gamma_{jik}\geq 0$, $\gamma_{kij}\geq 0$ are used.
Then we have $\det \mathcal{J}>0$ on $\mathcal{U}^H_{ijk}$.

By the connectivity of $\Omega^H_{ijk}$ and the continuity of the eigenvalues of $\Lambda^H_{ijk}$, we just need to prove
$\mathcal{J}$ is positive definite for some radius vector in $\Omega^H_{ijk}$.
By Corollary \ref{hyperbolic triangle inequality is valid for large s},
for sufficient large $s$, the radius vector $(s,s,s)\in \Omega^H_{ijk}$.
We shall prove $\mathcal{J}$ is positive definite for some $s$ large enough.
At $(s,s,s)$, we have
\begin{equation*}
\begin{aligned}
\mathcal{J}
=&\sinh^2s\cosh s
  \left(
   \begin{array}{ccc}
     \sinh l_i & 0 & 0 \\
     0 & \sinh l_j & 0 \\
     0 & 0 & \sinh l_k \\
   \end{array}
 \right)
 \left(
   \begin{array}{ccc}
     -1 & \cos \theta_k & \cos\theta_j \\
     \cos \theta_k & -1 & \cos \theta_i \\
     \cos \theta_j & \cos \theta_i & -1 \\
   \end{array}
 \right)\\
 &\times
  \left(
   \begin{array}{ccc}
     \frac{1}{\sinh l_i} & 0 & 0 \\
     0 & \frac{1}{\sinh l_j} & 0 \\
     0 & 0 & \frac{1}{\sinh l_k} \\
   \end{array}
   \right)
   \left(
     \begin{array}{ccc}
       0 & 1+I_i & 1+I_i \\
       1+I_j & 0 & 1+I_j \\
       1+I_k & 1+I_k & 0 \\
     \end{array}
   \right).
\end{aligned}
\end{equation*}
Write the above equation as $\mathcal{J}=\sinh^2s\cosh s N$. Then
we just need to prove that the leading $1\times1$ and $2\times2$ minor of $N$ is positive for some $s$ large enough.

For the leading $1\times1$  minor, we have
\begin{equation}\label{positivity of 2*2 minor 1}
\begin{aligned}
N_{11}=&\frac{\sinh l_i\cos \theta_k}{\sinh l_j}(1+I_j)+ \frac{\sinh l_i\cos\theta_j}{\sinh l_k}(1+I_k)\\
=&\frac{1}{\sinh^2l_j\sinh^2l_k}[(1+I_j)(\cosh l_i\cosh l_j-\cosh l_k)(\cosh^2l_k-1)\\
 &+(1+I_k)(\cosh l_i\cosh l_k-\cosh l_j)(\cosh^2l_j-1)]\\
=&\frac{(1+I_j)(1+I_k)\sinh^4 s}{\sinh^2l_j\sinh^2l_k}[2(1+I_i)(1+I_j)(1+I_k)\sinh^4s\\
 &+(6+6I_i+3I_j+3I_k+3I_iI_j+3I_iI_k+2I_jI_k-I_j^2-I_k^2)\sinh^2s+4(1+I_i)].
\end{aligned}
\end{equation}
Note that, by Corollary \ref{hyperbolic triangle inequality is valid for large s}, under the condition
$$2\sinh^2 s(1+I_i)(1+I_j)(1+I_j)>I_i^2+I_j^2+I_k^2-3,$$
the triangle inequalities are satisfied,
which implies
\begin{equation*}
\begin{aligned}
&\frac{\sinh l_i\cos \theta_k}{\sinh l_j}(1+I_j)+ \frac{\sinh l_i\cos\theta_j}{\sinh l_k}(1+I_k)\\
\geq&\frac{(1+I_j)(1+I_k)\sinh^4 s}{\sinh^2l_j\sinh^2l_k}\\
&\times[(3+6I_i+3I_j+3I_k+3I_iI_j+3I_iI_k+2I_jI_k+I_i^2)\sinh^2s+4(1+I_i)]\\
=&\frac{(1+I_j)(1+I_k)\sinh^4 s}{\sinh^2l_j\sinh^2l_k}[\left((1+I_i)(3+I_i)+2\gamma_{ijk}+3\gamma_{jik}+3\gamma_{kij}\right)\sinh^2s+4(1+I_i)].
\end{aligned}
\end{equation*}
Therefor the leading $1\times1$  minor of $N$ is positive by
the condition $I_{i}, I_j, I_k\in (-1, +\infty)$ and $\gamma_{ijk}\geq 0$, $\gamma_{jik}\geq 0$, $\gamma_{kij}\geq 0$.

Similar to (\ref{positivity of 2*2 minor 1}), we have
\begin{equation}\label{positivity of 2*2 minor 2}
\begin{aligned}
N_{22}=&\frac{\sinh l_j\cos \theta_k}{\sinh l_i}(1+I_i)+\frac{\sinh l_j\cos \theta_i}{\sinh l_k}(1+I_k)\\
=&\frac{(1+I_i)(1+I_k)\sinh^4 s}{\sinh^2l_i\sinh^2l_k}[2(1+I_i)(1+I_j)(1+I_k)\sinh^4s\\
 &+(6+3I_i+6I_j+3I_k+3I_iI_j+2I_iI_k+3I_jI_k-I_i^2-I_k^2)\sinh^2s+4(1+I_j)].
\end{aligned}
\end{equation}
Note that
\begin{equation}\label{positivity of 2*2 minor 3}
\begin{aligned}
N_{12}N_{21}
=&[-(1+I_i)+\frac{\sinh l_i\cos\theta_j}{\sinh l_k}(1+I_k)][-(1+I_j)+\frac{\sinh l_j\cos \theta_i}{\sinh l_k}(1+I_k)]\\
=&\frac{1}{\sinh^4l_k}[(1+I_k)\sinh l_k\sinh l_i\cos\theta_j-(1+I_i)\sinh^2l_k]\\
 &\times[(1+I_k)\sinh l_k\sinh l_j\cos\theta_i-(1+I_j)\sinh^2l_k]\\
=&\frac{1}{\sinh^4l_k}[(1+I_k)(\cosh l_i\cosh l_k-\cosh l_j)-(1+I_i)\sinh^2l_k]\\
 &\times[(1+I_k)(\cosh l_j\cosh l_k-\cosh l_i)-(1+I_j)\sinh^2l_k]\\
=&\frac{(1+I_k)^4\sinh^4s}{\sinh^4l_k}(1+I_i+I_j-I_k)^2,
\end{aligned}
\end{equation}
where $\cosh l_i=\cosh^2s+I_i\sinh^2s=1+(1+I_i)\sinh^2s$ is used in the last line.

Combining (\ref{positivity of 2*2 minor 1}), (\ref{positivity of 2*2 minor 2}), (\ref{positivity of 2*2 minor 3}), we have
the leading $2\times2$  minor of $N$ is
\begin{equation*}
\begin{aligned}
&\frac{(1+I_i)(1+I_j)(1+I_k)^2\sinh^8s}{\sinh^2l_i\sinh^2l_j\sinh^4l_k}\\
 &\times\big[2(1+I_i)(1+I_j)(1+I_k)\sinh^4s\\
 &\ \ \ \ +(6+6I_i+3I_j+3I_k+3I_iI_j+3I_iI_k+2I_jI_k-I_j^2-I_k^2)\sinh^2s+4(1+I_i)\big]\\
 &\times[2(1+I_i)(1+I_j)(1+I_k)\sinh^4s\\
 &\ \ \ \ +(6+3I_i+6I_j+3I_k+3I_iI_j+2I_iI_k+3I_jI_k-I_i^2-I_k^2)\sinh^2s+4(1+I_j)]\\
 &-\frac{(1+I_k)^4\sinh^4s}{\sinh^4l_k}(1+I_i+I_j-I_k)^2\\
=&\frac{(1+I_k)^2\sinh^4s}{\sinh^2l_i\sinh^2l_j\sinh^4l_k}\\
 &\times\{(1+I_i)(1+I_j)\sinh^4s[4(1+I_i)^2(1+I_j)^2(1+I_k)^2\sinh^8s\\
 &\ \ \ \ \ \ \ \ \ \ \ \ \ \ \ \ \ \ \ \ \ \ \ \ \ \ \ \ \ \ \ \ \ \ \ \ +A\sinh^6s+B\sinh^4s+C\sinh^2s+D]\\
 &\ \ \ \ -(1+I_k)^2(1+I_i+I_j-I_k)^2\sinh^2l_i\sinh^2l_j\},
\end{aligned}
\end{equation*}
where $A,B,C,D$ are polynomials of $I_i, I_j, I_k$.
Note that
$\sinh^2l_i=\cosh^2l_i-1=(1+I_i)\sinh^2s[2+(1+I_i)\sinh^2s],$
we have the leading $2\times2$  minor of $N$ is
\begin{equation*}
\begin{aligned}
&\frac{(1+I_i)(1+I_j)(1+I_k)^2\sinh^8s}{\sinh^2l_i\sinh^2l_j\sinh^4l_k}\\
 &\times\{4(1+I_i)^2(1+I_j)^2(1+I_k)^2\sinh^8s+A\sinh^6s+B\sinh^4s+C\sinh^2s+D\\
 &\ \ \ \ -(1+I_k)^2(1+I_i+I_j-I_k)^2[2+(1+I_i)\sinh^2s][2+(1+I_j)\sinh^2s]\}.
\end{aligned}
\end{equation*}
The term in the last two lines is a polynomial in $\sinh s$ with positive leading coefficient $4(1+I_i)^2(1+I_j)^2(1+I_k)^2$,
so for $s$ large enough, the leading $2\times2$  minor of $N$ is positive.

Combining with the fact that the determinant of $\mathcal{J}$ is positive, we have
the matrix $\Lambda^H_{ijk}$ is negative definite. This completes the proof.
\qed

\begin{remark}
The matrix $\mathcal{J}$ in (\ref{matrix J}) is the same matrix $M$ in the proof of Lemma 12
of Guo \cite{Guo}, where $M$ was proved to be positive definite for nonnegative inversive distance.
Here we produces another proof of the fact.
\end{remark}

\begin{remark}
If $I_{i}, I_j, I_k\in (-1, 1]$ and $\gamma_{ijk}\geq 0$, $\gamma_{jik}\geq 0$, $\gamma_{kij}\geq 0$,
the negative definiteness of $\Lambda^H_{ijk}$ was proved by Zhou \cite{Z} using the same method
as that of  Lemma \ref{hyperbolic negative definite}.
In this case, the negative definiteness of $\Lambda^H_{ijk}$ could be proved alternatively.
In fact, by direct but tedious calculations, we have
\begin{equation}\label{derivative of area}
\begin{aligned}
&\frac{\partial \theta_i}{\partial u_i}+\frac{\partial \theta_j}{\partial u_i}+\frac{\partial \theta_k}{\partial u_i}\\
=&\frac{1}{A(\cosh l_j+1)(\cosh l_k+1)}\cdot\\
&\{C_iS_i^2S_j^2(I_k^2-1)+S_i^2S_j^2C_k(I_k^2-1)-S_k(C_jS_i^2S_j\gamma_{ijk}+C_iS_iS_j^2\gamma_{jik})\\
&+C_kS_k[-S_i^2S_j(2C_iC_j+1)\gamma_{ijk}-(C_i^2+S_i^2)S_iS_j^2\gamma_{jik}]\\
&+S_k^2[-2C_iS_i^2S_j^2(I_iI_jI_k+1)-S_iS_j\gamma_{kij}(C_jS_i^2+C_i)+C_iS_i^2(I_j^2-1)+C_jS_i^2(I_j^2-1)]\}.
\end{aligned}
\end{equation}
In general, $\frac{\partial \theta_i}{\partial u_i}+\frac{\partial \theta_j}{\partial u_i}+\frac{\partial \theta_k}{\partial u_i}$
have no sign. However, if $I_{i}, I_j, I_k\in (-1, 1]$ and $\gamma_{ijk}\geq 0$, $\gamma_{jik}\geq 0$, $\gamma_{kij}\geq 0$,
we have
$\frac{\partial \theta_i}{\partial u_i}+\frac{\partial \theta_j}{\partial u_i}+\frac{\partial \theta_k}{\partial u_i}<0$
by (\ref{derivative of area}). Combining with Remark
\ref{positivity of hyperbolic derivative}, this implies $-\Lambda^H_{ijk}$ is diagonal dominant and then
$\Lambda^H_{ijk}$ is negative definite.
\end{remark}

Set
$$\Lambda^H=\frac{\partial (K_1,\cdots, K_N)}{\partial(u_1,\cdots,u_N)}=-\sum_{\triangle ijk\in F}\Lambda^H_{ijk},$$
where $\Lambda^H_{ijk}$ is extended by zeros to a $N\times N$
matrix so that $\Lambda^H_{ijk}$ acts on a vector $(v_1, \cdots, v_N)$ only on
the coordinates corresponding to vertices $v_i, v_j$ and $v_k$ in the triangle $\triangle ijk$.
Lemma \ref{hyperbolic symmetry lemma} and Lemma \ref{hyperbolic negative definite} have the following direct corollary.
\begin{corollary}\label{positivity of Lambda H}
Given a triangulated surface $(M, \mathcal{T}, I)$ with inversive distance $I>-1$ and
$\gamma_{ijk}\geq 0, \gamma_{jik}\geq 0, \gamma_{kij}\geq 0$
for any topological triangle $\triangle ijk\in F$.
Then the matrix $\Lambda^H=\frac{\partial (K_1,\cdots, K_N)}{\partial(u_1,\cdots,u_N)}$
is symmetric and positive definite
on $\mathcal{U}^H:=\cap_{\triangle ijk\in T}\mathcal{U}^H_{ijk}$ for the hyperbolic background geometry.
\end{corollary}
Guo \cite{Guo} once obtained a result paralleling to Corollary \ref{positivity of Lambda H} for $I\geq 0$.

By Lemma \ref{hyperbolic simply connect} and Lemma \ref{hyperbolic symmetry lemma}, we can define an energy function
$$\mathcal{E}_{ijk}(u)=\int_{u_0}^u\theta_idu_i+\theta_jdu_j+\theta_kdu_k$$
on $\mathcal{U}^H_{ijk}=\ln (\Omega^H_{ijk})$. Lemma \ref{hyperbolic negative definite} ensures that
$\mathcal{E}_{ijk}$ is locally concave on $\mathcal{U}^H_{ijk}$.
Define the Ricci potential as
\begin{equation}\label{hyperbolic total energy}
\mathcal{E}(u)=-\sum_{\triangle ijk\in T}\mathcal{E}_{ijk}(u)+\int_{u_0}^u\sum_{i=1}^N(2\pi-\overline{K}_i)du_i,
\end{equation}
then $\nabla_u\mathcal{E}=K-\overline{K}$ and $\mathcal{E}(u)$ is locally convex on $\mathcal{U}^H=\cap_{\triangle ijk\in T}\mathcal{U}^H_{ijk}$.
The local convexity of $\mathcal{E}$ implies the infinitesimal rigidity of $K$ with respect to $u$, which is the infinitesimal rigidity of hyperbolic inversive
distance circle packings.

\subsection{Global rigidity of hyperbolic inversive distance circle packings}

In this subsection, we shall prove the global rigidity of hyperbolic inversive distance circle packings under the condition $I\in (-1, +\infty)$ and
$\gamma_{ijk}\geq 0, \gamma_{jik}\geq 0, \gamma_{kij}\geq 0$
for any triangle $\triangle ijk\in F$.

By Corollary \ref{hyperbolic triangle extension lemma}, the functions $\theta_i, \theta_j, \theta_k$ defined on $\mathcal{U}^H_{ijk}$
could be continuously extended by constants to $\widetilde{\theta}_i, \widetilde{\theta}_j, \widetilde{\theta}_k$ defined on $\mathbb{R}^3$.
Using Theorem \ref{Luo's convex extention}, we have the following extension.

\begin{lemma}\label{hyperbolic extension lemma}
In the hyperbolic background geometry, for any triangle $\triangle ijk\in F$ with $I_i, I_j, I_k>-1$ and
$\gamma_{ijk}\geq 0$, $\gamma_{jik}\geq 0$, $\gamma_{kij}\geq 0$,
the function $\mathcal{E}_{ijk}(u)$ defined on $\mathcal{U}^H_{ijk}$
could be extended to the following function
\begin{equation}\label{extension function Euclid}
\widetilde{\mathcal{E}}_{ijk}(u)=\int_{u_0}^u\widetilde{\theta}_idu_i+\widetilde{\theta}_jdu_j+\widetilde{\theta}_kdu_k,
\end{equation}
which is a $C^1$-smooth concave function defined on $\mathbb{R}^3$ with
$$\nabla_u \widetilde{\mathcal{E}}_{ijk}=(\widetilde{\theta}_i, \widetilde{\theta}_j, \widetilde{\theta}_k)^T.$$
\end{lemma}

Using Lemma \ref{hyperbolic extension lemma}, we can prove the following global rigidity of hyperbolic inversive distance circle packing metrics,
which is the hyperbolic part of Theorem \ref{main theorem global rigidity for Curvature}.
.
\begin{theorem}\label{hyperbolic rigidity}
Given a triangulated surface $(M, \mathcal{T})$ with inversive distance $I\in (-1, +\infty)$ and
$
\gamma_{ijk}\geq 0, \gamma_{jik}\geq 0, \gamma_{kij}\geq 0
$
for any topological triangle $\triangle ijk\in F$. Then for any $\overline{K}\in C(V)$,
there is at most one hyperbolic inversive distance circle packing metric $r$ with $K(r)=\overline{K}$.
\end{theorem}

\proof
The Ricci energy function $\mathcal{E}(u)$ in (\ref{hyperbolic total energy}) could be extended from $\mathcal{U}^H$ to the whole space $\mathbb{R}^{N}$, where
$\mathcal{U}^H$ is the image of $\Omega^H$ under the map $u_i=\ln \tanh\frac{r_i}{2}$.
In fact, the function $\mathcal{E}_{ijk}(u)$ defined on $\mathcal{U}^H_{ijk}$
could be extended to $\widetilde{\mathcal{E}}_{ijk}(u)$ defined by (\ref{extension function Euclid})
on $\mathbb{R}^N$ by Lemma \ref{hyperbolic extension lemma}
and the second term $\int_{u_0}^u\sum_{i=1}^N(2\pi-\overline{K}_i)du_i$ in (\ref{hyperbolic total energy})
can be naturally defined on $\mathbb{R}^N$, then we have the following extension $\widetilde{\mathcal{E}}(u)$ defined on $\mathbb{R}^N$
of the Ricci potential function $\mathcal{E}(u)$
$$\widetilde{\mathcal{E}}(u)=-\sum_{\triangle ijk\in F}\widetilde{\mathcal{E}}_{ijk}(u)+\int_{u_0}^u\sum_{i=1}^N(2\pi-\overline{K}_i)du_i.$$
As $\widetilde{\mathcal{E}}_{ijk}(u)$ is $C^1$-smooth concave by Lemma \ref{hyperbolic extension lemma}
and $\int_{u_0}^u\sum_{i=1}^N(2\pi-\overline{K}_i)du_i$
is a well-defined convex function on $\mathbb{R}^N$,
we have $\widetilde{\mathcal{E}}(u)$ is a $C^1$-smooth convex function on $\mathbb{R}^N$. Furthermore,
\begin{equation*}
\nabla_{u_i}\widetilde{F}=-\sum_{\triangle ijk\in F}\widetilde{\theta}_i+2\pi-\overline{K}_i=\widetilde{K}_i-\overline{K}_i,
\end{equation*}
where $\widetilde{K}_i=2\pi-\sum_{\triangle ijk\in F}\widetilde{\theta}_i$.

If there are two different inversive distance circle packing metrics $\overline{r}_{A}, \overline{r}_{B}\in \Omega^H$ with the same combinatorial
Curvature $\overline{K}$, then
$\overline{u}_A=\ln \tanh\frac{\overline{r}_{A}}{2}\in \mathcal{U}^H$, $\overline{u}_B=\ln \tanh\frac{\overline{r}_{B}}{2}\in \mathcal{U}^H$
are both critical points of the extended Ricci potential $\widetilde{\mathcal{E}}(u)$.
It follows that
$$\nabla \widetilde{\mathcal{E}}(\overline{u}_A)=\nabla \widetilde{\mathcal{E}}(\overline{u}_B)=0.$$
Set
\begin{equation*}
\begin{aligned}
f(t)=&\widetilde{\mathcal{E}}((1-t)\overline{u}_A+t\overline{u}_B)\\
=&\sum_{\triangle ijk\in F}f_{ijk}(t)+\int_{u_0}^{(1-t)\overline{u}_A+t\overline{u}_B}\sum_{i=1}^N(2\pi-\overline{K}_i)du_i,
\end{aligned}
\end{equation*}
where
$$f_{ijk}(t)=-\widetilde{\mathcal{E}}_{ijk}((1-t)\overline{u}_A+t\overline{u}_B).$$
Then $f(t)$ is a $C^1$ convex function on $[0, 1]$ and $f'(0)=f'(1)=0$, which implies $f'(t)\equiv 0$ on $[0, 1]$.
Note that $\overline{u}_A$ belongs to the open set $\mathcal{U}^H$,
there exists $\epsilon>0$ such that $(1-t)\overline{u}_A+t\overline{u}_B\in \mathcal{U}^H$ for $t\in [0, \epsilon]$.
So $f(t)$ is smooth on $[0, \epsilon]$.

Note that $f(t)$ is $C^1$ convex on $[0, 1]$ and smooth on $[0, \epsilon]$.
$f'(t)\equiv 0$ on $[0, 1]$ implies that $f''(t)\equiv 0$ on $[0, \epsilon]$.
Note that, for $t\in [0, \epsilon]$,
\begin{equation*}
\begin{aligned}
f''(t)=(\overline{u}_A-\overline{u}_B) \Lambda^H  (\overline{u}_A-\overline{u}_B)^T,
\end{aligned}
\end{equation*}
where $\Lambda^H=-\sum_{\triangle ijk\in F}\Lambda^H_{ijk}$.
By Corollary \ref{positivity of Lambda H}, we have $\Lambda^H$ is positive definite and then $\overline{u}_A-\overline{u}_B=0$, which
implies that $\overline{r}_A=\overline{r}_B$. So there exists at most one hyperbolic inversive distance circle packing metric
with combinatorial curvature $\overline{K}$.
\qed

\subsection{Rigidity of combinatorial $\alpha$-curvature in hyperbolic background geometry}

We have the following global rigidity for $\alpha$-curvature with respect to hyperbolic inversive distance circle packing metrics
for inversive distance in $(-1, +\infty)$, which is the hyperbolic part of Theorem \ref{main theorem global rigidity for alpha curvature}.

\begin{theorem}\label{hyperbolic alpha rigidity}
Given a closed triangulated surface $(M, \mathcal{T})$ with inversive distance $I>-1$ and
$\gamma_{ijk}\geq 0, \gamma_{jik}\geq 0, \gamma_{kij}\geq 0$
for any topological triangle $\triangle ijk\in F$, $\overline{R}$
is a given function defined on the vertices  of  $(M, \mathcal{T})$.
If $\alpha\overline{R}\leq 0$,
there exists at most one hyperbolic inversive distance circle packing metric $\overline{r}\in \Omega^H$
with combinatorial $\alpha$-curvature $\overline{R}$.
\end{theorem}

As the proof of Theorem \ref{hyperbolic alpha rigidity} is almost parallel to that of Theorem \ref{hyperbolic rigidity} using the energy function
$$\widetilde{\mathcal{E}}_{\alpha}(u)=-\sum_{\triangle ijk\in F}\widetilde{\mathcal{E}}_{ijk}(u)+\int_{u_0}^u\sum_{i=1}^N(2\pi-\overline{R}_i\tanh^\alpha \frac{r_i}{2})du_i,$$
we omit the details of the proof here. Theorem \ref{hyperbolic alpha rigidity} is an generalization of Theorem \ref{hyperbolic rigidity}.
Specially, if $\alpha=0$, Theorem \ref{hyperbolic alpha rigidity} is reduced to Theorem \ref{hyperbolic rigidity}.\\

\textbf{Acknowledgements}\\[8pt]
The research of the author is supported by Hubei Provincial Natural Science Foundation of China under grant No. 2017CFB681,
Fundamental Research Funds for the Central Universities under grant No. 2042018kf0246 and
National Natural Science Foundation of China under grant No. 61772379 and No. 11301402.
Part of this work was done during the visit of Department of Mathematical Sciences, Tsinghua University.
The author thanks Professor Daguang Chen for his invitation and thanks Professor Feng Luo, Professor Xianfeng Gu, Professor Huabin Ge,
Professor Ze Zhou, Dr. Wai Yeung Lam and Xiang Zhu for communications
on related topics. The author thanks the referees for their careful reading of the manuscript and for their nice suggestions
which greatly improves the exposition of the paper.

(Xu Xu) School of Mathematics and Statistics, Wuhan University, Wuhan 430072, P.R. China

E-mail: xuxu2@whu.edu.cn\\[2pt]

\end{document}